\documentclass[11pt,a4paper]{amsart}
\usepackage{graphicx,url}

\theoremstyle{plain}
\newtheorem{theorem}{Theorem}[section]
\newtheorem{lemma}[theorem]{Lemma}
\newtheorem{corollary}[theorem]{Corollary}

\newtheorem{proposition}[theorem]{Proposition}
\theoremstyle{definition}

\newcommand{\thmlabel}[1]{\label{thm:#1}}
\newcommand{\thmref}[1]{Theorem~\ref{thm:#1}}
\newcommand{\lemlabel}[1]{\label{lem:#1}}
\newcommand{\lemref}[1]{Lemma~\ref{lem:#1}}

\newcommand{\threelemref}[3]{Lemmas~\ref{lem:#1}, \ref{lem:#2} and \ref{lem:#3}}
\newcommand{\eqnlabel}[1]{\label{eqn:#1}}
\newcommand{\eqnref}[1]{\eqref{eqn:#1}}
\newcommand{\Eqnref}[1]{Equation~\eqref{eqn:#1}}
\newcommand{\figlabel}[1]{\label{fig:#1}}
\newcommand{\figref}[1]{Figure~\ref{fig:#1}}

\newcommand{\seclabel}[1]{\label{sec:#1}}
\newcommand{\secref}[1]{Section~\ref{sec:#1}}
\newcommand{\twosecref}[2]{Sections~\ref{sec:#1} and \ref{sec:#2}}
\newcommand{\corlabel}[1]{\label{cor:#1}}
\newcommand{\corref}[1]{Corollary~\ref{cor:#1}}
\newcommand{\proplabel}[1]{\label{prop:#1}}
\newcommand{\propref}[1]{Proposition~\ref{prop:#1}}
\newcommand{\twopropref}[2]{Propositions~\ref{prop:#1} and \ref{prop:#2}}

\newcommand{\V}[1]{\ensuremath{\protect|#1|}}
\newcommand{\E}[1]{\ensuremath{\protect\|#1\|}}
\newcommand{\CR}[1]{\ensuremath{\protect\textsf{\textup{cr}}(#1)}}
\newcommand{\RCR}[1]{\ensuremath{\protect\overline{\textsf{\textup{cr}}}(#1)}}
\newcommand{\CCR}[1]{\ensuremath{\protect\textsf{\textup{cr}}^{\star}(#1)}}
\newcommand{\FFF}{\ensuremath{\mathcal{F}}}
\newcommand{\PPP}{\ensuremath{\mathcal{P}}}
\newcommand{\SSS}[1][\gamma]{\ensuremath{\mathbb{S}_{#1}}}
\newcommand{\tpw}[1]{\ensuremath{\textup{\textsf{tpw}}(#1)}}
\newcommand{\tw}[1]{\ensuremath{\textup{\textsf{tw}}(#1)}}

\date{April 21, 2006, Revised: \today}
\title[Planar Decompositions and the Crossing Number]{Planar Decompositions and the Crossing Number of Graphs with an Excluded Minor}

\author{David R. Wood}
\address{\newline
Departament de Matem{\'a}tica Aplicada II\newline 
Universitat Polit{\`e}cnica de Catalunya\newline 
Barcelona, Spain}
\email{\tt david.wood@upc.es}

\author{Jan Arne Telle}
\address{\newline
Department of Informatics\newline
The University of Bergen\newline 
Bergen, Norway}
\email{\tt Jan.Arne.Telle@ii.uib.no}

\thanks{This paper was published in the \emph{New York J. Math.} 13:117--146, 2007 (\url{http://nyjm.albany.edu/j/2007/13-8.html}). This version incorporates the journal version plus \twopropref{Strong3DecompK5}{StrongOmegaDecompK5}. An extended abstract of this paper was published in the \emph{Proceedings of the 14th International Symposium on Graph Drawing} (GD '06), \emph{Lecture Notes in Computer Science},
vol.~4372, pp.~150--161, Springer, 2007.}
 	
\thanks{The research of David Wood is supported by a Marie Curie Fellowship of the European Community under contract 023865, and by the projects 
MEC MTM2006-01267 and DURSI 2005SGR00692.}



\newcommand{\Figure}[4][htb]{
\begin{figure}[#1]
	\begin{center}#3\end{center}
	\caption{\figlabel{#2}#4}
\end{figure}
}

\newcommand{\ceil}[1]{\ensuremath{\protect\lceil#1\rceil}}

\newcommand{\floor}[1]{\ensuremath{\protect\lfloor#1\rfloor}}

\newcommand{\half}{\ensuremath{\protect\tfrac{1}{2}}}

\newcommand{\third}{\ensuremath{\protect\tfrac{1}{3}}}

\newcommand{\twothirds}{\ensuremath{\protect\tfrac{2}{3}}}

\begin{document}

\begin{abstract}
Tree decompositions of graphs are of fundamental importance in structural and algorithmic graph theory. Planar decompositions generalise tree decompositions by allowing an arbitrary planar graph to index the decomposition. We prove that every graph that excludes a fixed graph as a minor has a planar decomposition with bounded width and a linear number of bags. 

The crossing number of a graph is the minimum number of crossings in a drawing of the graph in the plane. We prove that planar decompositions are intimately related to the crossing number. In particular, a graph with bounded degree has linear crossing number if and only if it has a planar decomposition with bounded width and a linear number of bags. It follows from the above result about planar decompositions that every graph with bounded degree and an excluded minor has linear crossing number. 

Analogous results are proved for the convex and rectilinear crossing numbers. In particular, every graph with bounded degree and bounded tree-width has linear convex crossing number, and every $K_{3,3}$-minor-free graph with bounded degree has linear rectilinear crossing number.
\end{abstract}

\maketitle
\newpage
%\tableofcontents
%\newpage
\section{Introduction}
\label{sec:Intro}
%%%%%%%%%%%%%%%%%%%%%%%%%%%%%%%%%%%%%%%%%%%%%%%%%%%%%%%%%%%%%%%%%%%%%%%%%%%%%%%%

The \emph{crossing number} of a graph\footnote{We consider graphs $G$ that are undirected, simple, and finite. Let $V(G)$ and $E(G)$ respectively be the vertex and edge sets of $G$. Let $\V{G}:=|V(G)|$ and $\E{G}:=|E(G)|$. For each vertex $v$ of a graph $G$, let $N_G(v):=\{w\in V(G):vw\in E(G)\}$ be the neighbourhood of $v$ in $G$. The \emph{degree} of $v$ is $|N_G(v)|$. Let $\Delta(G)$ be the maximum degree of a vertex of $G$.} $G$, denoted by \CR{G}, is the minimum number of crossings in a drawing\footnote{A \emph{drawing} of a graph represents each vertex by a distinct point in the plane, and represents each edge by a simple closed curve between its endpoints, such that the only vertices an edge intersects are its own endpoints, and no three edges intersect at a common point (except at a common endpoint). A \emph{crossing} is a point of intersection between two edges (other than a common endpoint). A drawing with no crossings is \emph{plane}. A graph is \emph{planar} if it has a plane drawing.} of $G$ in the plane; see \cite{EG-AMM73, PachToth-JCTB00, Vrto} for surveys. The crossing number is an important measure of the non-planarity of a graph \cite{Szekely-DM04}, with applications in discrete and computational geometry \cite{PachSharir-CPC98, Szekely-CPC97} and VLSI circuit design \cite{BL84, Leighton83, Leighton84}. In information visualisation, one of the most important measures of the quality of a graph drawing is the number of crossings \cite{Purchase-JVLC98, Purchase97, PCJ97}. 

Computing the crossing number is $\mathcal{NP}$-hard \cite{GJ-SJDM83}, and remains so for simple cubic graphs \cite{Hliney-JCTB06, PSS05}. Moreover, the exact or even asymptotic crossing number is not known for specific graph families, such as complete graphs \cite{RT-AMM97}, complete bipartite graphs \cite{Nahas-EJC03, RS-JGT96, RT-AMM97}, and cartesian products \cite{AR-JCTB04, Bokal-JCTB07, GS-JGT04, RT-DCG95}. 

Given that the crossing number seems so difficult, it is natural to focus on asymptotic bounds rather than exact values. The `crossing lemma', conjectured by Erd{\H{o}}s and Guy~\cite{EG-AMM73} and first proved \footnote{A remarkably simple probabilistic proof of the crossing lemma was found by Chazelle, Sharir and Welzl~\cite{Proofs3}. See \cite{Montaron-JGT05, PRTT-DCG06} for recent improvements.} by Leighton~\cite{Leighton83} and Ajtai~et~al.~\cite{Ajtai82}, gives such a lower bound. It states that 
\begin{equation*}
\CR{G}\geq \E{G}^3/\,64\V{G}^2
\end{equation*}
for every graph $G$ with $\E{G}\geq4\V{G}$. Other general lower bound techniques that arose out of the work of Leighton~\cite{Leighton83, Leighton84} include the bisection/cutwidth method \cite{DV-JGAA03, PSS-Algo96, SS-CPC94, SSSV-Algo96} and the embedding method \cite{SSSV-AM96, SS-CPC94}. 

Upper bounds on the crossing number of general families of graphs have been less studied, and are the focus of this paper. Obviously $\CR{G}\leq\binom{\E{G}}{2}$ for every graph $G$. A family of graphs has \emph{linear}\footnote{If the crossing number of a graph is linear in the number of edges then it is also linear in the number of vertices. To see this, let $G$ be a graph with $n$ vertices and $m$ edges. Suppose that $\CR{G}\leq cm$. If $m<4n$ then $\CR{G}\leq4cn$ and we are done.  Otherwise  $\CR{G}\geq m^3/64n^2$ by the crossing lemma. Thus $m\leq8\sqrt{c}n$ and $\CR{G}\leq8c^{3/2}n$.} crossing number if 
\begin{equation*}
\CR{G}\leq c\,\V{G}
\end{equation*}
for every graph $G$ in the family, for some constant $c$. For example, Pach and T\'{o}th~\cite{PachToth-GD05} proved that graphs of bounded genus\footnote{Let \SSS\ be the orientable surface with $\gamma\geq0$ handles. An \emph{embedding} of a graph in \SSS\ is a crossing-free drawing in \SSS. A \emph{$2$-cell embedding} is an embedding in which each region of the surface (bounded by edges of the graph) is an open disk. The (\emph{orientable}) \emph{genus} of a graph $G$ is the minimum $\gamma$ such that $G$ has a $2$-cell embedding in \SSS. In what follows, by a \emph{face} we mean the set of vertices on the boundary of the face. Let $F(G)$ be the set of faces in an embedded graph $G$. See the monograph by Mohar and Thomassen~\cite{MoharThom} for a thorough treatment of graphs on surfaces.} and bounded degree have linear crossing number. Our main result states that bounded-degree graphs that exclude a fixed graph as a minor\footnote{Let $vw$ be an edge of a graph $G$. Let $G'$ be the graph obtained by identifying the vertices $v$ and $w$, deleting loops, and replacing parallel edges by a single edge. Then $G'$ is obtained from $G$ by \emph{contracting} $vw$. A graph $H$ is a \emph{minor} of a graph $G$ if $H$ can be obtained from a subgraph of $G$ by contracting edges. A family of graphs \FFF\ is \emph{minor-closed} if $G\in\FFF$ implies that every minor of $G$ is in \FFF. \FFF\ is \emph{proper} if it is not the family of all graphs. A deep theorem of Robertson and Seymour~\cite{RS-GraphMinorsXX-JCTB04} states that every proper minor-closed family can be characterised by a finite family of excluded minors. Every proper minor-closed family is a subset of the $H$-minor-free graphs for some graph $H$. We thus focus on minor-closed families with one excluded minor.} have linear crossing number.

\begin{theorem}
\thmlabel{CrossingMinorFree}
For every graph $H$ and integer $\Delta$, there is a constant $c=c(H,\Delta)$, such that every $H$-minor-free graph $G$ with maximum degree at most $\Delta$ has crossing number $\CR{G}\leq c\,\V{G}$.
\end{theorem}

\thmref{CrossingMinorFree} implies the above-mentioned result of Pach and T\'{o}th~\cite{PachToth-GD05}, since graphs of bounded genus exclude a fixed graph as a minor (although the dependence on $\Delta$ is different in the two proofs; see \secref{Surfaces}). Moreover, \thmref{CrossingMinorFree} is stronger than the above-mentioned result of Pach and T\'{o}th~\cite{PachToth-GD05}, since there are graphs with a fixed excluded minor and unbounded genus\footnote{Since the genus of a graph equals the sum of the genera of its biconnected components, it is trivial to construct $1$-connected graphs that exclude a fixed minor, yet have unbounded genus. There are highly connected examples as well: For fixed $p$, the complete bipartite graph $K_{p,n}$ is $p$-connected, has no $K_{p+2}$-minor, yet has unbounded genus \cite{Ringel65}. There are examples with bounded degree as well. Seese and Wessel~\cite{SW-JCTB89} constructed a family of graphs, each with no $K_8$-minor and maximum degree $5$, and with unbounded genus.}. For other recent work on minors and crossing number see \cite{BCSV-ENDM, BFM-SJDM06, BFW, GS-JGT01, GRS-EJC04, Hliney-JCTB03, Hliney-JCTB06, Negami-JGT01, PSS05}.

Note that the assumption of bounded degree in \thmref{CrossingMinorFree} is unavoidable. For example, $K_{3,n}$ has no $K_5$-minor, yet has $\Omega(n^2)$ crossing number \cite{RS-JGT96, Nahas-EJC03}. Conversely, bounded degree does not by itself guarantee linear crossing number. For example, a random cubic graph on $n$ vertices has $\Omega(n)$ bisection width \cite{CE-BAMS89, DDSW-TCS03}, which implies that it has $\Omega(n^2)$ crossing number \cite{DV-JGAA03, Leighton83}. 

The proof of \thmref{CrossingMinorFree} is based on \emph{planar decompositions}, which are introduced in \twosecref{Decompositions}{Manipulating}. This combinatorial structure generalises tree decompositions by allowing an arbitrary planar graph to index the decomposition. We prove that planar decompositions and the crossing number are intimately related (\secref{Key}). In particular, a graph with bounded degree has linear crossing number if and only if it has a planar decomposition with bounded width and linear order (\thmref{CrossingNumberCharacterisation}). We study planar decompositions of: $K_5$-minor-free graphs (\secref{K5}), graphs embedded in surfaces (\secref{Surfaces}), and finally graphs with an excluded minor (\secref{Minors}). One of the main contributions of this paper is to prove that every graph that excludes a fixed graph as a minor has a planar decomposition with bounded width and linear order. \thmref{CrossingMinorFree} easily follows.

\subsection{Complementary Results}
%%%%%%%%%%%%%%%%%%%%%%%%%%%%%%%%%%%%%%%%%%%%%%%%%%%%%%%%%%%%%%%%%%%%%%%%%%%

A graph drawing is \emph{rectilinear} (or \emph{geometric}) if each edge is represented by a straight line-segment. The \emph{rectilinear crossing number} of a graph $G$, denoted by \RCR{G}, is the minimum number of crossings in a rectilinear drawing of $G$; see \cite{AAK-Computing06, BD-JGT93, BDG-DM03, LVWW04, RT-AMM97, SW-AMM94}. A rectilinear drawing is \emph{convex} if the vertices are positioned on a circle. The \emph{convex} (or \emph{outerplanar}, \emph{circular}, or \emph{1-page book}) \emph{crossing number} of a graph $G$, denoted by \CCR{G}, is the minimum number of crossings in a convex drawing of $G$; see \cite{CSSV-EJC04, Riskin-BICA03, SSSV04}. Obviously 
\begin{equation*}
\CR{G}\leq\RCR{G}\leq\CCR{G}
\end{equation*}
for every graph $G$. \emph{Linear} rectilinear and \emph{linear} convex crossing numbers are defined in an analogous way to linear crossing number.

It is unknown whether an analogue of \thmref{CrossingMinorFree} holds for rectilinear crossing number\footnote{The crossing number and rectilinear crossing number are not related in general. In particular, for every integer $k\geq4$, Bienstock and Dean~\cite{BD-JGT93} constructed a graph $G_k$ with crossing number $4$ and rectilinear crossing number $k$. It is easily seen that $G_k$ has no $K_{14}$-minor. However, the maximum degree of $G_k$ increases with $k$. Thus $G_k$ is not a counterexample to an analogue of \thmref{CrossingMinorFree} for rectilinear crossing number.}. On the other hand, we prove that $K_{3,3}$-minor-free graphs with bounded degree have linear rectilinear crossing number (\secref{K33}). 

\begin{theorem}
\thmlabel{RectilinearCrossingK33}
For every integer $\Delta$, there is a constant $c=c(\Delta)$, such that every $K_{3,3}$-minor-free graph $G$ with maximum degree at most $\Delta$ has rectilinear crossing number $\RCR{G}\leq c\,\V{G}$.
\end{theorem}

An analogue of \thmref{CrossingMinorFree} for convex crossing number does not hold, even for planar graphs, since Shahrokhi~et~al.~\cite{SSSV04} proved that the $n\times n$ planar grid $G_n$ (which has maximum degree $4$) has convex crossing number $\Omega(\V{G_n}\log \V{G_n})$. It is natural to ask which property of the planar grid forces up the convex crossing number. In some sense, we show that tree-width\footnote{Tree-width is a minor-monotone parameter that is defined in \secref{Decompositions}.} is one answer to this question. In particular, $G_n$ has tree-width $n$. More generally, we prove that every graph with large tree-width has many crossings on some edge in every convex drawing (\propref{ConvexTreewidth}). On the other hand, we prove that graphs with bounded tree-width and bounded degree have linear convex crossing number (\secref{Partitions}). 

\begin{theorem}
\thmlabel{ConvexCrossingTreewidth}
For all integers $k$ and $\Delta$, there is a constant $c=c(k,\Delta)$, such that every graph $G$ with tree-width at most $k$ and maximum degree at most $\Delta$ has convex crossing number $\RCR{G}\leq c\,\V{G}$.
\end{theorem}

Again, the assumption of bounded degree in \thmref{ConvexCrossingTreewidth} is unavoidable since $K_{3,n}$ has tree-width $3$.

%Example: take $K_n$ and replace each vertex by a cycle. We get a 3-regular graph with $n(n-1)$ vertices, which I expect has $\Omega(n^4)$ crossing number.

%%%%%%%%%%%%%%%%%%%%%%%%%%%%%%%%%%%%%%%%%%%%%%%%%%%%%%%%%%%%%%%%%%%%%%%%%%%%%%%%
\section{Graph Decompositions}
\seclabel{Decompositions}
%%%%%%%%%%%%%%%%%%%%%%%%%%%%%%%%%%%%%%%%%%%%%%%%%%%%%%%%%%%%%%%%%%%%%%%%%%%%%%%%

%For $S\subseteq V(G)$, let $D(S)$ be the subgraph of $D$ induced by the bags that contain $S$. Let $D(v):=D(\{v\})$ for each vertex $v$ of $G$, and $D(vw):=D(\{v,w\})$ for each edge $vw$ of $G$. 

Let $G$ and $D$ be graphs, such that each vertex of $D$ is a set of vertices of $G$ (called a \emph{bag}). Note that we allow distinct vertices of $D$ to be the same set of vertices in $G$; that is, $V(D)$ is a multiset. For each vertex $v$ of $G$, let $D(v)$ be the subgraph of $D$ induced by the bags that contain $v$. Then $D$ is a \emph{decomposition} of $G$ if:
\begin{itemize}
\item $D(v)$ is connected and nonempty for each vertex $v$ of $G$, and 
\item $D(v)$ and $D(w)$ touch\footnote{Let $A$ and $B$ be subgraphs of a graph $G$. Then $A$ and $B$ \emph{intersect} if $V(A)\cap V(B)\ne\emptyset$, and $A$ and $B$ \emph{touch} if they intersect or $v\in V(A)$ and $w\in V(B)$ for some edge $vw$ of $G$. } for each edge $vw$ of $G$.
\end{itemize}

Decompositions, when $D$ is a tree, were introduced by Robertson and Seymour~\cite{RS-GraphMinorsII-JAlg86}. Diestel and K{\"u}hn~\cite{DiestelKuhn-DAM05}\footnote{A decomposition was called a \emph{connected decomposition} by Diestel and K{\"u}hn~\cite{DiestelKuhn-DAM05}. Similar definitions were introduced by Agnew~\cite{Agnew05}. Decompositions can also be defined in terms of the lexicographic product. For a graph $G$ and integer $k\geq 1$, the lexicographic product $G\cdot K_k$ is the graph with vertex set $V(G)\times[k]$, where $(v,i)(w,j)$ is an edge of $G\cdot K_k$ if and only if $vw\in E(G)$, or $v=w$ and $i\ne j$. That is, each vertex of $G$ is `blown up' by a copy of $K_k$, and each edge of $G$ is `blown up' by a copy of $K_{k,k}$. It is easily seen that $D$ is a decomposition of $G$ with width $k$ if and only if $G$ is a minor of $D\cdot K_k$. With this viewpoint, a similar result to \lemref{QuadraticDecomp} was obtained by Ne{\v{s}}et{\v{r}}il and Ossona De Mendez~\cite{NesOdM-GradI}, who observed that every graph is a minor of $D\cdot K_2$ for some planar graph $D$.} first generalised the definition for arbitrary graphs $D$.

Let $D$ be a decomposition of a graph $G$. The \emph{width} of $D$ is the maximum cardinality of a bag. The number of bags that contain a vertex $v$ of $G$ is the \emph{spread} of $v$ in $D$. The \emph{spread} of $D$ is the maximum spread of a vertex of $G$. The \emph{order} of $D$ is the number of bags. $D$ has \emph{linear order} if $\V{D}\leq c\,\V{G}$ for some constant $c$. If the graph $D$ is a tree, then the decomposition $D$ is a \emph{tree decomposition}. If the graph $D$ is a cycle, then the decomposition $D$ is a \emph{cycle decomposition}. The decomposition $D$ is \emph{planar} if the graph $D$ is planar. The \emph{genus} of the decomposition $D$ is the genus of the graph $D$.

Note that decompositions generalise minors since $D$ is a decomposition of $G$ with width $1$ if and only if $G$ is a minor of $D$. 

A decomposition $D$ of a graph $G$ is \emph{strong} if $D(v)$ and $D(w)$ intersect for each edge $vw$ of $G$. The \emph{tree-width} of $G$, denoted by \tw{G}, is $1$ less than the minimum width of a strong tree decomposition of $G$. For example, a graph has tree-width $1$ if and only if it is a forest. Graphs with tree-width $2$ (called \emph{series-parallel}) are planar, and are characterised as those graphs with no $K_4$-minor. Tree-width is particularly important in structural and algorithmic graph theory; see the surveys \cite{Bodlaender-TCS98, Reed-AlgoTreeWidth03}. 

For applications to crossing number, tree decompositions are not powerful enough: even the $n\times n$ planar grid has tree-width $n$. We show in \secref{Key} that planar decompositions are the right type of decomposition for applications to crossing number. It is tempting to define the `planar-width' of a graph $G$ to be the minimum width in a planar decomposition of $G$. However, by the following lemma of Diestel and K{\"u}hn~\cite{DiestelKuhn-DAM05}, every graph would then have bounded planar-width. We include the proof for completeness.

\begin{lemma}[\cite{DiestelKuhn-DAM05}]
\lemlabel{QuadraticDecomp}
Every graph $G$ has a strong planar decomposition of width $2$, spread $\V{G}+1$, and order $\binom{\V{G}+1}{2}$. 
\end{lemma}

\begin{proof}
Let $n:=\V{G}$ and say $V(G)=\{1,2,\dots,n\}$. Define a graph $D$ with vertex set $V(D):=\{\{i,j\}:1\leq i\leq j\leq n\}$ and edge set $E(D):=\{\{i,j\}\{i+1,j\}:1\leq i\leq n-1,1\leq j\leq n\}$. Then $D$ is a planar subgraph of the $n\times n$ grid; see \figref{GridDecomposition}. For each vertex $i$ of $G$, the set of bags that contain $i$ is $\{\{i,j\}:1\leq j\leq n\}$, which induces a (connected) $n$-vertex path in $D$. For each edge $ij$ of $G$, the bag $\{i,j\}$ contains $i$ and $j$. Therefore $D$ is a strong decomposition of $G$. The width is $2$, since each bag has two vertices. Each vertex is in $n+1$ bags.
\end{proof}

%\begin{floatingfigure}{55mm}\includegraphics{GridDecomposition}\caption{A strong planar decomposition of $K_6$ with width $2$ and order $21$; the subgraph $D(3)$ is highlighted.}\end{floatingfigure} 

\Figure{GridDecomposition}{\includegraphics{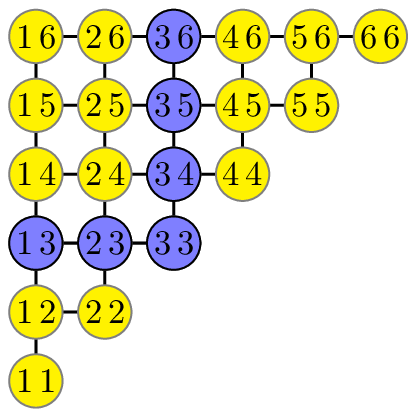}}{A strong planar decomposition of $K_6$ with width $2$ and order $21$; the subgraph $D(3)$ is highlighted.}

%\comment{Note that we can delete the $\{i,i\}$ bag, and connect each $\{i-1,i\}$ to $\{i,i+1\}$, and still get the result.}

The planar decomposition in \lemref{QuadraticDecomp} has large order (quadratic in \V{G}). The remainder of this paper focuses on planar decompositions with linear order. 

Strong tree decompositions are the most widely studied decompositions in the literature. This paper focuses on decompositions that are not necessarily strong. One advantage is that every graph obviously has a decomposition isomorphic to itself (with width $1$). On the other hand, if $G$ has a strong decomposition $D$ of width $k$, then 
\begin{equation}
\eqnlabel{NumberEdgesStrongDecomp}
\E{G}\leq\tbinom{k}{2}\,\V{D}\enspace.
\end{equation}
It follows that if $G$ has a strong decomposition isomorphic to itself then the width is at least $\sqrt{\frac{2\E{G}}{\V{G}}}$, which is unbounded for dense graphs, as observed by Diestel and K{\"u}hn~\cite{DiestelKuhn-DAM05}. Note that if $G$ has a (non-strong) decomposition $D$ of width $k$, then 
\begin{equation}
\eqnlabel{NumberEdgesDecomp}
\E{G}\leq k^2\,\E{D}+\tbinom{k}{2}\,\V{D}\enspace.
\end{equation}

Every tree $T$ satisfies the Helly property: every collection of pairwise intersecting subtrees of $T$ have a vertex in common. It follows that if a tree $T$ is a strong decomposition of $G$ then every clique\footnote{A \emph{clique} of a graph $G$ is a set of pairwise adjacent vertices in $G$. The maximum cardinality of a clique of $G$ is denoted by $\omega(G)$.} of $G$ is contained in some bag of $T$. Other graphs do not have this property. It will be desirable (for performing $k$-sums in \secref{Manipulating}) that (non-tree) decompositions have a similar property. We therefore introduce the following definitions. 

For $p\geq0$, a \emph{$p$-clique} is a clique of cardinality $p$. A \emph{$(\leq p)$-clique} is a clique of cardinality at most $p$. For $p\geq2$, a decomposition $D$ of a graph $G$ is a  \emph{$p$-decomposition} if each $(\leq p)$-clique of $G$ is a subset of some bag of $D$, or is a subset of the union of two adjacent bags of $D$. An $\omega(G)$-decomposition of $G$ is called an \emph{$\omega$-decomposition}. A $p$-decomposition $D$ of $G$ is \emph{strong} if each $(\leq p)$-clique of $G$ is a subset of some bag of $D$. Observe that a (strong) $2$-decomposition is the same as a (strong) decomposition, and a (strong) $p$-decomposition also is a (strong) $q$-decomposition for all $q\in[2,p]$.

\section{Manipulating Decompositions}
\seclabel{Manipulating}
%%%%%%%%%%%%%%%%%%%%%%%%%%%%%%%%%%%%%%%%%%%%%%%%%%

In this section we describe four tools for manipulating graph decompositions that are repeatedly used in the remainder of the paper. 

\subsection*{Tool \#$1$. Contracting a Decomposition:}
%%%%%%%%%%%%%%%%%%%%%%%%%%%%%%%%%%%%%%%%%%%%%%%%%%%%%%%%%

Our first tool describes the effect of contracting an edge in a decomposition.

\begin{lemma}[Contraction]
\lemlabel{ContractDecomp}
Suppose that $D$ is a planar (strong) $p$-decomposition of a graph $G$ with width $k$. Say $XY$ is an edge of $D$. Then the decomposition $D'$ obtained by contracting the edge $XY$ into the vertex $X\cup Y$ is a planar (strong) $p$-decomposition of $G$ with width $\max\{k,|X\cup Y|\}$. In particular, if 
$|X\cup Y|\leq k$ then $D'$ also has width $k$. 
\end{lemma}

\begin{proof}
Contracting edges preserves planarity. Thus $D'$ is planar. Contracting edges preserves connectiveness. Thus $D'(v)$ is connected for each vertex $v$ of $G$. Contracting the edge $XY$ obviously maintains the required properties for each $(\leq p)$-clique of $G$. 
\end{proof}

\lemref{ContractDecomp} can be used to decrease the order of a decomposition at the expense of increasing the width. The following observation is a corollary of \lemref{ContractDecomp}.

\begin{corollary}
\corlabel{ContractMatching}
Suppose that $D$ is a (strong) $p$-decomposition of a graph $G$ with width $k$, and that $D$ has a matching\footnote{A \emph{matching} is a set of pairwise disjoint edges.} $M$. The decomposition obtained from $D$ by contracting $M$  is a (strong) $p$-decomposition of $G$ with width at most $2k$ and order $\V{D}-|M|$.\qed
\end{corollary}

\begin{lemma}
\lemlabel{ReduceOrder}
Suppose that a graph $G$ has a (strong) planar $p$-decomposition $D$ of width $k$ and order at most $c\V{G}$ for some $c\geq1$. Then $G$ has a (strong) planar $p$-decomposition of width $c'k$ and order $\V{G}$, for some $c'$ depending only on $c$. 
\end{lemma}

\begin{proof}
Without loss of generality, $D$ is a planar triangulation. Biedl~et~al.~\cite{BDDFK-DM04} proved that every planar triangulation on $n$ vertices has a matching of at least $\frac{n}{3}$ edges. Applying this result to $D$, and by  \corref{ContractMatching}, $G$ has a (strong) planar $p$-decomposition of width at most $2k$ and order at most $\twothirds\V{D}$. By induction, for every integer $i\geq1$, $G$ has a (strong) planar $p$-decomposition of width $2^ik$ and order at most $\max\{(\twothirds)^i\V{D},1\}$. With $i:=\ceil{\log_{3/2}c}$, the assumption that $\V{D}=c\V{G}$ implies that $G$ has a (strong) planar $p$-decomposition of width $2^ik$ and order $\V{G}$.
\end{proof}

\subsection*{Tool \#$2$. Composing Decompositions:}
%%%%%%%%%%%%%%%%%%%%%%%%%%%%%%%%%%%%%%%%%%%%%%%%%%%%%%%%%%%%%

Our second tool describes how two decompositions can be composed.

\begin{lemma}[Composition Lemma]
\lemlabel{Composition}
Suppose that $D$ is a (strong) $p$-decomposition of a graph $G$ with width $k$, and that $J$ is a decomposition of $D$ with width $\ell$. Then $G$ has a (strong) $p$-decomposition isomorphic to $J$ with width $k\ell$. 
\end{lemma}

\begin{proof} 
Let $J'$ be the graph isomorphic to $J$ that is obtained by renaming each bag $Y\in V(J)$ by $Y':=\{v\in V(G):v\in X\in Y\text{ for some }X\in V(D)\}$. There are at most $\ell$ vertices $X\in Y$, and at most $k$ vertices $v\in X$. Thus each bag of $J'$ has at most $k\ell$ vertices. 

First we prove that $J'(v)$ is connected for each vertex $v$ of $G$. Let $A'$ and $B'$ be two bags of $J'$ that contain $v$. Let $A$ and $B$ be the corresponding bags in $D$. Thus $v\in X_1$ and $v\in X_t$ for some bags $X_1,X_t\in V(D)$ such that $X_1\in A$ and $X_t\in B$ (by the construction of $J'$). Since $D(v)$ is connected, there is a path $X_1,X_2,\dots,X_t$ in $D$ such that $v$ is in each $X_i$. In particular, each $X_iX_{i+1}$ is an edge of $D$. Now $J(X_i)$ and $J(X_{i+1})$ touch in $J$. Thus there is path in $J$ between any vertex of $J$ that contains $X_1$ and any vertex of $J$ that contains $X_t$, such that every bag in the path contains some $X_i$. In particular, there is a path $P$ in $J$ between $A$ and $B$ such that every bag in $P$ contains some $X_i$. Let $P':=\{Y':Y\in P\}$. Then $v\in Y'$ for each bag $Y'$ of $P'$ (by the construction of $J'$). Thus $P'$ is a connected subgraph of $J'$ that includes $A'$ and $B'$, and $v$ is in every such bag. Therefore $J'(v)$ is connected. 

%Since $J$ is a decomposition of $D$, for $i\in[1,t-1]$, there is bag $Y_i$ of $J$ that contains $X_i$ and $X_{i+1}$. Thus for $i\in[2,t-1]$, $X_i$ is in $Y_{i-1}$ and $Y_i$, and there is a path $P_i$ in $J$ between $Y_{i-1}$ and $Y_i$, such that each bag of $P_i$ contains $X_i$ (since $J(X_i)$ is connected). Let $P_i':=\{Y':Y\in P_i\}$. Then $v\in Y'$ for each bag $Y'$ of $P'_i$ (by the construction of $J'$). Thus the union of the $P'_i$ form a connected subgraph of $J'$ that includes $A'$ and $B'$, and $v$ is in every such bag. Therefore $J'(v)$ is connected. 

It remains to prove that for each $(\leq p)$-clique $C$ of $G$,\\
\hspace*{1.5em} (a) $C$ is a subset of some bag of $J'$, or\\ 
\hspace*{1.5em} (b) $C$ is a subset of the union of two adjacent bags of $J'$.\\ Moreover, we must prove that if $D$ is strong then case (a) always occurs.  Since $D$ is a $p$-decomposition of $G$, \\
\hspace*{1.5em} (1) $C\subseteq X$ for some bag $X\in V(D)$, or \\
\hspace*{1.5em} (2) $C\subseteq X_1\cup X_2$ for some edge $X_1X_2$ of $D$.\\
First suppose that case (1) applies, which always occurs if $D$ is strong. Then there is some bag $Y\in V(J)$ such that $X\in Y$ (since $X$ is a vertex of $D$ and $J$ is a decomposition of $D$). Thus $C\subseteq Y'$ by the construction of $J'$. Now suppose that case (2) applies. Then $D(X_1)$ and $D(X_2)$ touch in $D$. That is, $X_1$ and $X_2$ are in a common bag of $D$, or $X_1\in Y_1$ and $X_2\in Y_2$ for some edge $Y_1Y_2$ of $D$. If $X_1$ and $X_2$ are in a common bag $Y$, then since $C\subseteq X_1\cup X_2$, we have $C\subseteq Y'$ by the construction of $J'$; that is, case (a) occurs. Otherwise, $X_1\in Y_1$ and $X_2\in Y_2$ for some edge $Y_1Y_2$ of $D$. Then $C\cap X_1\subseteq Y_1'$ and $C\cap X_2\subseteq Y_2'$. Since $C\subseteq X_1\cup X_2$ we have $C\subseteq Y_1'\cup Y_2'$; that is, case (b) occurs. 
\end{proof}

\subsection*{Tool \#$3$. $\omega$-Decompositions:}
%%%%%%%%%%%%%%%%%%%%%%%%%%%%%%%%%%%%%%%%%%%%%%%%%%%%%%%%%

The third tool converts a decomposition into an $\omega$-decomposition with a small increase in the width. A graph $G$ is \emph{$d$-degenerate} if every subgraph of $G$ has a vertex of degree at most $d$. 

\begin{lemma}
\lemlabel{Degen}
Every $d$-degenerate graph $G$ has a strong $\omega$-decomposition isomorphic to $G$ of width at most $d+1$.
\end{lemma}

\begin{proof} 
It is well known (and easily proved) that $G$ has an acyclic orientation\footnote{If each edge of a graph $G$ is directed from one endpoint to the other, then we speak of an \emph{orientation} of $G$ with arc set $A(G)$. An orientation with no directed cycle is \emph{acyclic}. For each vertex $v$ of an orientation of $G$, let $N^-_G(v):=\{w\in V(G):(w,v)\in A(G)\}$ and $N^+_G(v):=\{w\in V(G):(v,w)\in A(G)\}$. The \emph{indegree} and \emph{outdegree} of $v$ are $|N^-_G(v)|$ and $N^+_G(v)|$ respectively. A \emph{sink} is a vertex with outdegree $0$.} such that each vertex has indegree at most $d$. Replace each vertex $v$ by the bag $\{v\}\cup N^-_G(v)$. Every subgraph of $G$ has a \emph{sink}. Thus every clique is a subset of some bag. The set of bags that contain a vertex $v$ are indexed by $\{v\}\cup N^+_G(v)$, which induces a connected subgraph in $G$. Thus we have a strong $\omega$-decomposition. Each bag has cardinality at most $d+1$.
\end{proof}

\begin{lemma}
\lemlabel{DegenDegen}
Suppose that $D$ is a decomposition of a $d$-degenerate graph $G$ of width $k$. Then $G$ has a strong $\omega$-decomposition isomorphic to $D$ of width $k(d+1)$. 
\end{lemma}

\begin{proof}
By \lemref{Degen}, $G$ has a strong $\omega$-decomposition isomorphic to $G$ of width $d$. By \lemref{Composition}, $G$ has a strong $\omega$-decomposition isomorphic to $D$ with width $k(d+1)$.
\end{proof}

In \lemref{DegenDegen}, the `blow-up' in the width is bounded by a constant factor for the graphs that we are interested in: even in the most general setting, $H$-minor-free graphs are $c\V{H}\sqrt{\log\V{H}}$-degenerate for some constant $c$ \cite{Kostochka82,Thomason84,Thomason01}.

\subsection*{Tool \#$4$. Clique-Sums of Decompositions:}
%%%%%%%%%%%%%%%%%%%%%%%%%%%%%%%%%%%%%%%%%%%%%%%%%%%%%%%%%

Our fourth tool describes how to determine a planar decomposition of a clique-sum of two graphs, given planar decompositions of the summands\footnote{Lea{\~n}os and Salazar~\cite{LS06} recently proved some related results on the additivity of crossing numbers.}. Let $G_1$ and $G_2$ be disjoint graphs. Suppose that $C_1$ and $C_2$ are $k$-cliques of $G_1$ and $G_2$ respectively, for some integer $k\geq0$. Let $C_1=\{v_1,v_2,\dots,v_k\}$ and $C_2=\{w_1,w_2,\dots,w_k\}$. Let $G$ be a graph obtained from $G_1\cup G_2$ by identifying $v_i$ and $w_i$ for each $i\in[1,k]$, and possibly deleting some of the edges $v_iv_j$. Then $G$ is a \emph{$k$-sum} of $G_1$ and $G_2$ \emph{joined} at $C_1=C_2$. An $\ell$-sum for some $\ell\leq k$ is called a \emph{$(\leq k)$-sum}. For example, if $G_1$ and $G_2$ are planar then it is easily seen that every $(\leq 2)$-sum of $G_1$ and $G_2$ is also planar, as illustrated in \figref{PlanarSum}.

\Figure{PlanarSum}{\includegraphics{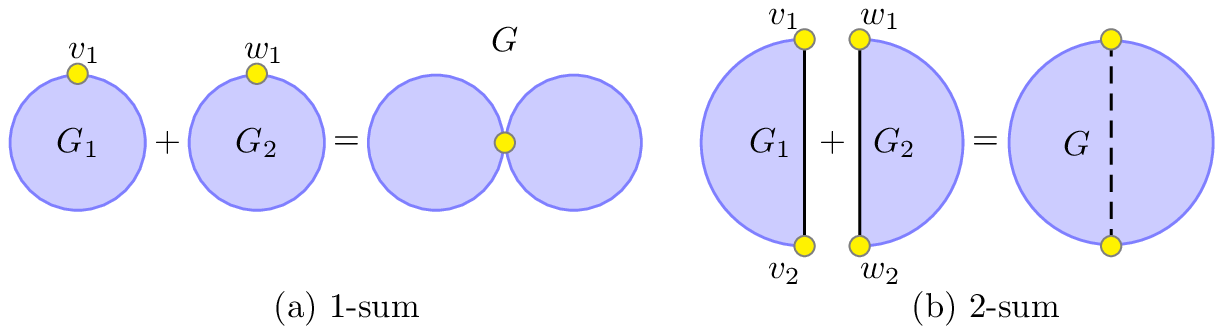}}{Clique-sums of planar graphs. In (b) we can assume that the edges $v_1v_2$ and $w_1w_2$ are respectively on the outerfaces of $G_1$ and $G_2$.}

\begin{lemma}[Clique Sums]
\lemlabel{Sum}
Suppose that for integers $p\leq q$, a graph $G$ is a $(\leq p)$-sum of graphs $G_1$ and $G_2$, and each $G_i$ has a (strong) planar $q$-decomposition $D_i$ of width $k_i$. Then $G$ has a (strong) planar $q$-decomposition of width $\max\{k_1,k_2\}$ and order $\V{D_1}+\V{D_2}$. 
\end{lemma}

\begin{proof}
Let $C:=V(G_1)\cap V(G_2)$. Then $C$ is a $(\leq p)$-clique, and thus a $(\leq q)$-clique, of both $G_1$ and $G_2$. Thus for each $i$,  \\
\hspace*{1.5em} (1) $C\subseteq X_i$ for some bag $X_i$ of $D_i$, or \\
\hspace*{1.5em} (2) $C\subseteq X_i\cup Y_i$ for some edge $X_iY_i$ of $D_i$.\\ If (1) is applicable, which is the case if $D_i$ is strong, then consider $Y_i:=X_i$ in what follows. 

Let $D$ be the graph obtained from the disjoint union of $D_1$ and $D_2$ by adding edges $X_1X_2$, $X_1Y_2$, $Y_1X_2$, and $Y_1Y_2$. By considering $X_1Y_1$ to be on the outerface of $G_1$ and $X_2Y_2$ to be on the outerface of $G_2$, observe that $D$ is planar, as illustrated in \figref{Sum}.

\Figure{Sum}{\includegraphics{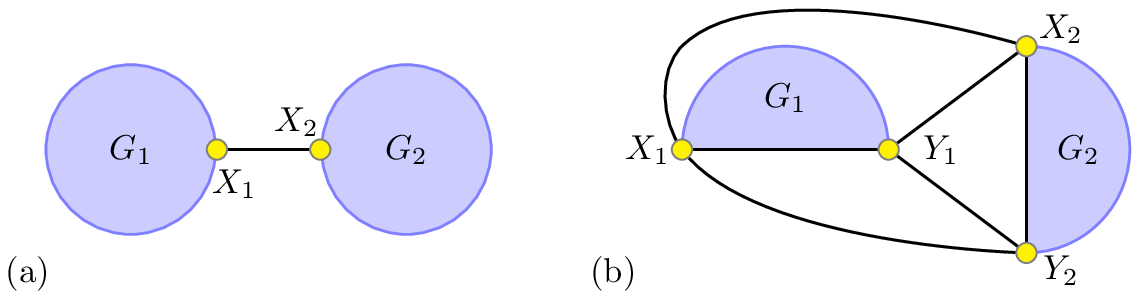}}{Sum of (a) strong planar decompositions, (b) planar decompositions.}

We now prove that $D(v)$ is connected for each vertex $v$ of $G$. If $v\not\in V(G_1)$ then $D(v)=D_2(v)$, which is connected. If $v\not\in V(G_2)$ then $D(v)=D_1(v)$, which is connected. Otherwise, $v\in C$. Thus $D(v)=D_1(v)\cup D_2(v)$. Since $v\in X_1\cup Y_1$ and $v\in X_2\cup Y_2$, and $X_1,Y_1,X_2,Y_2$ induce a connected subgraph ($\subseteq K_4$) in $D$, we have that $D(v)$ is connected. 

Each $(\leq q)$-clique $B$ of $G$ is a $(\leq q)$-clique of $G_1$ or $G_2$. Thus $B$ is a subset of some bag of $D$, or $B$ is a subset of the union of two adjacent bags of $D$. Moreover, if $D_1$ and $D_2$ are both strong, then $B$ is a subset of some bag of $D$. Therefore $D$ is a $q$-decomposition of $G$, and if $D_1$ and $D_2$ are both strong then $D$ is also strong. The width and order of $D$ are obviously as claimed. 
\end{proof}

A graph can be obtained by repeated $(\leq 1)$-sums of its biconnected
components. Thus \lemref{Sum} with $p=1$ implies:

\begin{corollary}
Let $G$ be a graph with biconnected components $G_1,G_2,\dots,G_t$. Suppose that each $G_i$ has a (strong) planar $q$-decomposition of width $k_i$ and order $n_i$. Then $G$ has a (strong) planar $q$-decomposition of width $\max_ik_i$ and order $\sum_in_i$.\qed
\end{corollary}

Note that in the proof of \lemref{Sum}, if  $X_1\subseteq X_2$ (for example) then, by \lemref{ContractDecomp}, we can contract the edge $X_1X_2$ in $D$ and merge the corresponding bags. The width is unchanged and the order is decreased by $1$. This idea is repeatedly used in the remainder of the paper. 

%%%%%%%%%%%%%%%%%%%%%%%%%%%%%%%%%%%%%%%%%%%%%%%%%%%%%%%%%%%%%%%%%%%%%%%%%%%%%%%%
\section{Planar Decompositions and the Crossing Number}
\seclabel{Key}
%%%%%%%%%%%%%%%%%%%%%%%%%%%%%%%%%%%%%%%%%%%%%%%%%%%%%%%%%%%%%%%%%%%%%%%%%%%%%%%%

The following lemma is the key link between planar decompositions and the crossing number of a graph.

\begin{lemma}
\lemlabel{DecompToCrossing}
Suppose that $D$ is a planar decomposition of a graph $G$ with width $k$. 
Then the crossing number of $G$ satisfies
\begin{equation*}
\CR{G}
\;\leq\;2\,\Delta(G)^2\!\!\sum_{X\in V(D)}\!\!\!\!\tbinom{|X|+1}{2}
\;\leq\;k(k+1)\,\Delta(G)^2\,\V{D}\enspace.
\end{equation*}
Moreover, if $s(v)$ is the spread of each vertex $v$ of $G$ in $D$, then $G$ has a drawing with the claimed number of crossings, in which each edge $vw$ is represented by a polyline with at most $s(v)+s(w)-2$ bends. 
\end{lemma}

\begin{proof}
By the F{\'a}ry-Wagner Theorem \cite{Fary48, Wagner36}, $D$ has a rectilinear drawing with no crossings. Let $R_\epsilon(X)$ be the open disc of radius $\epsilon>0$ centred at each vertex $X$ in the drawing of $D$. For each edge $XY$ of $D$, let $R_\epsilon(XY)$ be the union of all segments with one endpoint in $R_\epsilon(X)$ and one endpoint in $R_\epsilon(Y)$. For some $\epsilon>0$,
\begin{enumerate}
\item[(a)] $R_\epsilon(X)\cap R_\epsilon(Y)=\emptyset$ for all distinct bags $X$ and $Y$ of $D$, and 
\item[(b)] $R_\epsilon(XY)\cap R_\epsilon(AB)=\emptyset$ for all edges $XY$ and $AB$ of $D$ that have no endpoint in common.
\end{enumerate}

For each vertex $v$ of $G$, choose a bag $S_v$ of $D$ that contains $v$.
For each vertex $v$ of $G$, choose a point $p(v)\in R_\epsilon(S_v)$, and for each bag $X$ of $D$, choose a set $P(X)$ of $\sum_{v\in X}\deg_G(v)$ points in $R_\epsilon(X)$, so that if $\PPP=\cup\{P(X):X\in V(D)\}\cup\{p(v):v\in V(G)\}$ then:
\begin{enumerate}
\item[(c)] no two points in \PPP\ coincide, 
\item[(d)] no three points in \PPP\ are collinear, and 
\item[(e)] no three segments, each connecting two points in \PPP, cross at a common point.
\end{enumerate}
The set \PPP\ can be chosen iteratively since each disc $R_\epsilon(X)$ is $2$-dimensional\footnote{Let $Q$ be a nonempty set of points in the plane. Then $Q$ is \emph{$2$-dimensional} if it contains a disk of positive radius; $Q$ is \emph{$1$-dimensional} if it is not $2$-dimensional but contains a finite curve; otherwise $Q$ is \emph{$0$-dimensional}.}, but the set of excluded points is $1$-dimensional.

Draw each vertex $v$ at $p(v)$. For each edge $vw$ of $G$, a simple polyline  \begin{equation*}
L(vw)=(p(v),x_1,x_2,\dots,x_a,y_1,y_2,\dots,y_b,p(w))
\end{equation*}
(defined by its endpoints and bends) is a \emph{feasible} representation of $vw$ if:
\begin{enumerate}
\item $a\in[0,s(v)-1]$ and $b\in[0,s(w)-1]$,
\item each bend $x_i$ is in $P(X_i)$ for some bag $X_i$ containing $v$,
\item each bend $y_i$ is in $P(Y_i)$ for some bag $Y_i$ containing $w$, 
\item the bags $S_v,X_1,X_2,\dots,X_a,Y_1,Y_2,\dots,Y_b,S_w$ are distinct\\ (unless $S_v=S_w$ in which case $a=b=0$), and
\item consecutive bends in $L(vw)$ occur in adjacent bags of $D$.
\end{enumerate}
Since $D(v)$ and $D(w)$ touch, there is a feasible polyline that represents $vw$. 

A drawing of $G$ is \emph{feasible} if every edge of $G$ is represented by a feasible polyline, and no two bends coincide. Since each $|P(X)|=\sum_{v\in X}\deg(v)$, there is a feasible drawing. In particular, no edge passes through a vertex by properties (c)--(e), and no three edges have a common crossing point by property (e). By property (1), each edge $vw$ has at most $s(v)+s(w)-2$ bends.

Now choose a feasible drawing that minimises the total (Euclidean) length of the edges (with $\{p(v):v\in V(G)\}$ and $\{P(X):X\in V(D)\}$ fixed). 

%We claim that no two edges of $G$ with a common endpoint cross. Suppose on the contrary that edges $vu$ and $vw$ cross. The endpoints of the two segments of $vu$ and $vw$ that cross form a convex quadrilateral, of which the two segments are diagonals. Replace the diagonals with the appropriate pair of opposite sides of the quadrilateral so that two polylines are obtained, one between $p(v)$ to $p(u)$, and the other between $p(v)$ and $p(w)$; see \figref{RemoveCrossing}. We obtain a feasible drawing of $G$ with less total length. This contradiction proves that no two edges with a common endpoint cross. 

By properties (a), (b) and (2)--(5), each segment in a feasible drawing is 
contained within $R_\epsilon(X)$ for some bag $X$ of $D$, or within $R_\epsilon(XY)$ for some edge $XY$ of $D$. Consider a crossing in $G$ between edges $vw$ and $xy$. Since $D$ is drawn without crossings, the crossing point is contained within $R_\epsilon(X)$ for some bag $X$ of $D$, or within $R_\epsilon(XY)$ for some edge $XY$ of $D$. Thus some endpoint of $vw$, say $v$, and some endpoint of $xy$, say $x$, are in a common bag $X$. In this case, charge the crossing to the $5$-tuple $(vw,v,xy,x,X)$. Observe that the number of such $5$-tuples is 
\begin{equation*}
\sum_{X\in V(D)}\;\sum_{v,x\in X}\deg_G(v)\cdot\deg_G(x)\enspace.
\end{equation*}

At most four crossings are charged to each $5$-tuple $(vw,v,xy,x,X)$, since by property (4), each of $vw$ and $xy$ have at most two segments that intersect $R_\epsilon(X)$ (which might pairwise cross). We claim that, in fact, at most two crossings are charged to each such $5$-tuple. 

Suppose on the contrary that at least three crossings are charged to some $5$-tuple $(vw,v,xy,x,X)$. Then two segments of $vw$ intersect $R_\epsilon(X)$ and two segments of $xy$ intersect $R_\epsilon(X)$. In particular, $p(v)\not\in R_\epsilon(X)$ and $p(x)\not\in R_\epsilon(X)$, and $vw$ and $xy$ each have a bend in $R_\epsilon(X)$. Let $(r_1,r_2,r_3)$ be the $2$-segment polyline in the representation of $vw$, where $r_2$ is the bend of $vw$ in $R_\epsilon(X)$. Let $(t_1,t_2,t_3)$ be the $2$-segment polyline in the representation of $xy$, where $t_2$ is the bend of $xy$ in $R_\epsilon(X)$. Since at least three crossings are charged to $(vw,v,xy,x,X)$, in the set of segments $\{r_1r_2,r_2r_3,t_1t_2,t_2t_3\}$, at most one pair of segments, one from $vw$ and one from $xy$, do not cross. Without loss of generality, $t_1t_2$ and $r_2r_3$ are this pair. Observe that the crossing segments $r_1r_2$ and $t_1t_2$ are the diagonals of the convex quadrilateral $r_1t_2r_2t_1$. Replace the segments $r_1r_2$ and $t_1t_2$ by the segments $r_1t_2$ and $t_1r_2$, which are on opposite sides of the quadrilateral. Thus the combined length of $r_1t_2$ and $t_1r_2$ is less than the combined length of $r_1r_2$ and $t_1t_2$. Similarly, replace the segments $r_3r_2$ and $t_3t_2$ by the segments $r_3t_2$ and $t_3r_2$. We obtain a feasible drawing of $G$ with less total length. This contradiction proves that at most two crossings are charged to each $5$-tuple $(vw,v,xy,x,X)$. 

Thus the number of crossings is at most twice the number of $5$-tuples. Therefore the number of crossings is at most 
\begin{equation*}
2\sum_{X\in V(D)}\;\sum_{v,x\in X}\deg_G(v)\cdot\deg_G(x)
\;\leq\;
2\,\Delta(G)^2\,\sum_{X\in V(D)}\tbinom{|X|+1}{2}
\enspace.
\end{equation*}
\end{proof}

\Figure{RemoveCrossing}{\includegraphics{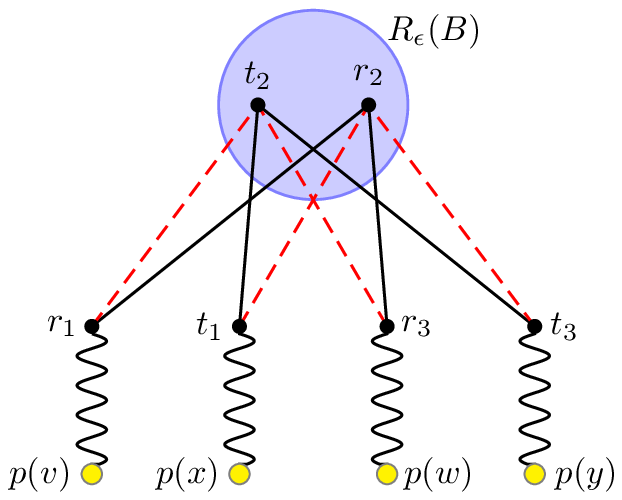}}{In the proof of \lemref{DecompToCrossing}, to shorten the total edge length the crossed segments are replaced by the dashed segments.}

Note that the bound on the crossing number in \lemref{DecompToCrossing} is within a constant factor of optimal for the complete graph. An easy generalisation of \lemref{QuadraticDecomp} proves that for all $n\geq k\geq 2$, $K_n$ has a strong planar decomposition of width $k$ and order at most $c(\frac{n}{k})^2$ for some constant $c$. Thus \lemref{DecompToCrossing} implies that  $\CR{K_n}\leq c k(k+1)\,\Delta(K_n)^2\,(\frac{n}{k})^2\leq cn^4$, which is within a constant factor of optimal \cite{RT-AMM97}.

%AAK-Computing06,RT97,Guy60}.

%%%%%%%%%%%%%%%%%%%%%%%%%%%%%%%%%%%%%%%%%

The following result is converse to \lemref{DecompToCrossing}.

\begin{lemma}
\lemlabel{CrossingToDecomp}
Let $G$ be a graph with $n$ isolated vertices. Suppose that $G$ has a drawing with $c$ crossings in which $q$ non-isolated vertices are not incident to a crossed edge. Then $G$ has a planar decomposition of width $2$ and order $\ceil{\frac{n}{2}}+q+c$, and $G$ has a strong planar decomposition of width $2$ and order $\ceil{\frac{n}{2}}+c+\E{G}$. 
\end{lemma}

\begin{proof}
First, pair the isolated vertices of $G$. Each pair can form one bag in a decomposition of width $2$, adding $\ceil{\frac{n}{2}}$ to the order. Now assume that $G$ has no isolated vertices. 

We first construct the (non-strong) decomposition. Arbitrarily orient each edge of $G$. Let $D$ be the planar graph obtained from the given drawing of $G$ by replacing each vertex $v$ by the bag $\{v\}$, and replacing each crossing between arcs $(v,w)$ and $(x,y)$ of $G$ by a degree-$4$ vertex $\{v,x\}$. Thus an arc $(v,w)$ of $G$ is replaced by some path $\{v\}\{v,x_1\}\{v,x_2\}\dots\{v,x_r\}\{w\}$ in $D$. In particular, $D(v)$ and $D(w)$ touch at the edge $\{v,x_r\}\{w\}$. Moreover, $D(v)$ is a (connected) tree for each vertex $v$ of $G$. Thus $D$ is a decomposition of $G$. Each bag contains at most two vertices. The order is $\V{G}+c$. For each vertex $v$ of $G$ that is incident to some crossed edge, $\{v\}\{v,x\}$ is an edge in $D$ for some vertex $x$. Contract the edge $\{v\}\{v,x\}$ in $D$ and merge the bags. By \lemref{ContractDecomp}, $D$ remains a planar decomposition of width $2$. The order is now $q+c$. 

Now we make $D$ strong. For each edge $vw$ of $G$, there is an edge $XY$ of $D$ where $D(v)$ and $D(w)$ touch. That is, $v\in X$ and $w\in Y$. Replace $XY$ by the path $X\{v,w\}Y$. Now each edge of $G$ is in some bag of $D$, and $D$ is strong. This operation introduces a further \E{G} bags. Thus the order is $q+c+\E{G}$. For each non-isolated vertex $v$ that is not incident to a crossed edge, choose an edge $vw$ incident to $v$. Then $\{v\}\{v,w\}$ is an edge in $D$; contract this edge and merge the bags. By \lemref{ContractDecomp}, $D$ remains a strong decomposition of $G$ with width $2$. The order is now $c+\E{G}$.
\end{proof}

The following two special cases of \lemref{CrossingToDecomp} are of particular interest.

\begin{corollary}
Every graph $G$ has a planar decomposition of width $2$ and order $\V{G}+\CR{G}$. Every graph $G$ has a strong planar decomposition of width $2$ and order $\V{G}+\E{G}+\CR{G}$. \qed
\end{corollary}

\begin{corollary}
Every planar graph $G$ with $n$ isolated vertices has a strong planar decomposition of width $2$ and order $\ceil{\frac{n}{2}}+\E{G}\leq 3\V{G}-6$. \qed
\end{corollary}

The relationship between planar decompositions and graphs with linear crossing number is summarised as follows. 

\begin{theorem}
\thmlabel{CrossingNumberCharacterisation}
The following are equivalent for a graph $G$ of bounded degree:
\begin{enumerate}
\item $\CR{G}\leq c_1\V{G}$ for some constant $c_1$,
\item $G$ has a planar decomposition with width $c_2$ and order $\V{G}$ for some constant $c_2$,
\item $G$ has a planar decomposition with width $2$ and order $c_3\V{G}$ for some constant $c_3$.
\end{enumerate}
\end{theorem}

\begin{proof}
\lemref{CrossingToDecomp} implies that (1) $\Rightarrow$ (3). 
\lemref{ReduceOrder} implies that (3) $\Rightarrow$ (2).
\lemref{DecompToCrossing} implies that (2) $\Rightarrow$ (1).
\end{proof}

%\twolemref{DecompToCrossing}{CrossingToDecomp} imply that the crossing number and planar decompositions are intimately related in the following sense.  

%\begin{corollary}\corlabel{CrossingDecomposition}Let \FFF\ be a family of graphs with bounded degree. Then \FFF\ has linear crossing number if and only if every graph in \FFF\ has a planar decomposition with bounded width and linear order. \qed\end{corollary}

%%%%%%%%%%%%%%%%%%%%%%%%%%%%%%%%%%%%%%%%%%%%%%%%%%%%%%%%%%%%%%%%%%%%%%%%%%%%%%%%
\section{$K_5$-Minor-Free Graphs}
\seclabel{K5}
%%%%%%%%%%%%%%%%%%%%%%%%%%%%%%%%%%%%%%%%%%%%%%%%%%%%%%%%%%%%%%%%%%%%%%%%%%%%%%%%

In this section we prove the following upper bound on the crossing number.

\begin{theorem}
\thmlabel{CrossingsK5}
Every $K_5$-minor-free graph $G$ has crossing number 
\begin{equation*}
\CR{G}\;<\;\tfrac{20}{3}\,\Delta(G)^2\,\V{G}\enspace.
\end{equation*}
\end{theorem}

The proof of \thmref{CrossingsK5} is based on \thmref{PlanarOmegaDecompK5} below, in which we construct $\omega$-decompositions of $K_5$-minor-free graphs $G$. Since $\omega(G)\leq4$ and each clique can be spread over two bags, it is natural to consider $\omega$-decompositions of width $2$.

\begin{theorem}
\thmlabel{PlanarOmegaDecompK5}
Every $K_5$-minor-free graph $G$ with $\V{G}\geq3$ has a planar $\omega$-decomposition of width $2$ and order at most $\frac{4}{3}\V{G}-2$, with at most $\V{G}-2$ bags of cardinality $2$, and at most $\frac{1}{3}\V{G}$ bags of cardinality $1$. 
\end{theorem}

The proof of \thmref{PlanarOmegaDecompK5} is based on the following classical theorem of Wagner~\cite{Wagner37} and the two following lemmas. Here $V_8$ is the the $8$-vertex {M}\"obius ladder, which is the graph obtained from the $8$-cycle by adding an edge between each pair of antipodal vertices; see \figref{W}(a). 

\begin{theorem}[\cite{Wagner37}]
\thmlabel{CharacterisationK5}
A graph is $K_5$-minor-free if and only if it can be obtained from planar graphs and $V_8$ by $(\leq 3)$-sums. 
\end{theorem}

\Figure{W}{\includegraphics{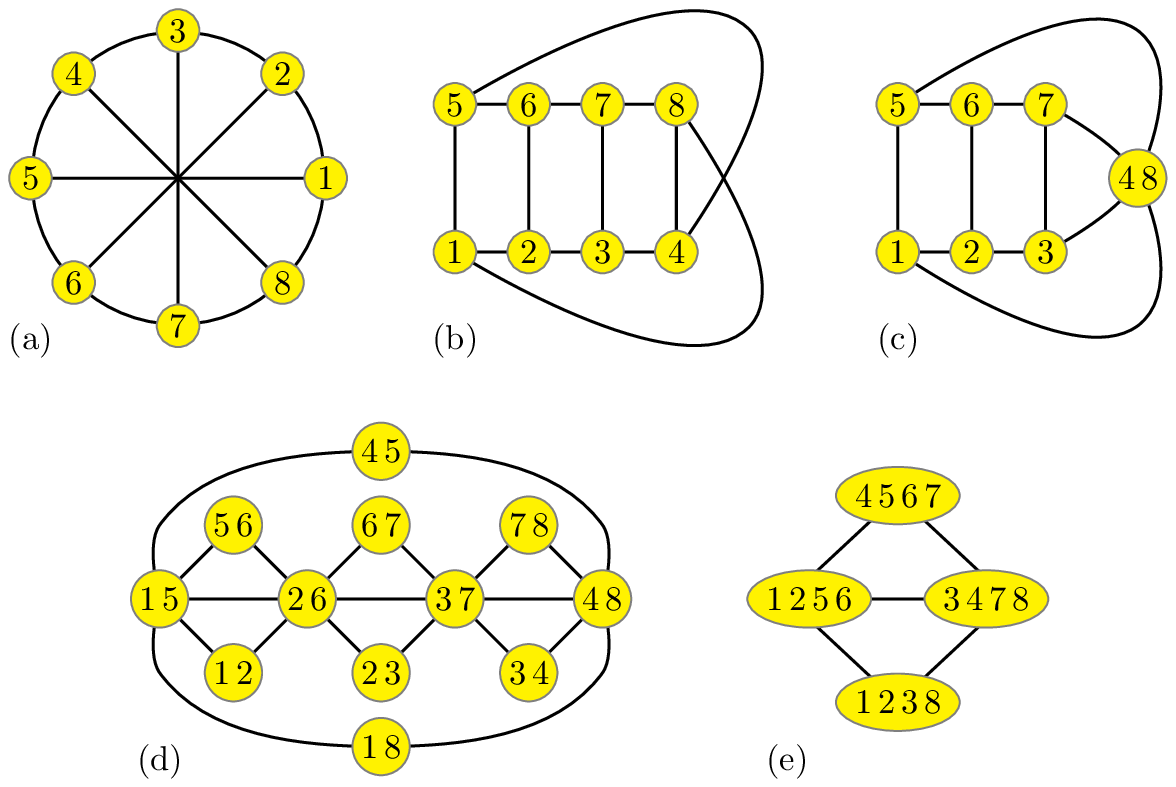}}{(a) The graph $V_8$. (b) Drawing of $V_8$ with one crossing.\newline (c) Planar $\omega$-decomposition of $V_8$ with width $2$ and order $7$.\newline (d) Strong planar decomposition of $V_8$ with width $2$ and order $12$.\newline (e)  Strong planar $\omega$-decomposition of $V_8$ with width $4$ and order $4$.}

%%%%%%%%%%%%%%%%%%%%%%%%%%%%%%%%%%%%%%

\begin{lemma}
\lemlabel{EdgePartitionK5}
Every $K_5$-minor-free graph $G$ with $\V{G}\geq3$ has a partition of $E(G)$ into three sets $E^1,E^2,E^3$ such that:
\begin{itemize}
\item each of $E^1,E^2,E^3$ has at most $\V{G}-2$ edges, 
\item every triangle has one edge in each of $E^1,E^2,E^3$,
\item if a subgraph $H$ of $G$ is isomorphic to $V_8$, then $E(H)\cap E^j$ is a perfect matching in $H$ for all $j$.
\end{itemize}
Moreover, if $G$ is edge-maximal (with no $K_5$-minor), then every vertex is incident to an edge in $E^j$ and an edge in $E^\ell$ for some $j\ne \ell$.
\end{lemma}

\begin{proof}
By \thmref{CharacterisationK5}, we need only consider the following three cases.

\textbf{Case (a).} $G$ is planar: Let $G'$ be a planar triangulation of $G$. By the Four-Colour Theorem \cite{RSST97}, $G'$ has a proper vertex-colouring with colours $a,b,c,d$. Now determine a Tait edge-colouring \cite{Tait1880a}. 
Let $E^1$ be the set of edges of $G'$ whose endpoints are coloured $ab$ or $cd$. Let $E^2$ be the set of edges of $G'$ whose endpoints are coloured $ac$ or $bd$. Let $E^3$ be the set of edges of $G'$ whose endpoints are coloured $ad$ or $bc$. Since the vertices of each triangle are $3$-coloured, the edges of each triangle are in distinct $E^j$. In particular, the edges of each face of $G'$ are in distinct $E^j$. Each edge of $G'$ is in two of the $2\V{G}-4$ faces of $G'$. Thus $|E^j|=\V{G}-2$. The sets $E^j\cap E(G)$ thus satisfy the first two properties for $G$. Since $V_8$ is nonplanar, $G$ has no $V_8$ subgraph, and the third property is satisfied vacuously. Finally, if $G$ is edge-maximal, then $G'=G$, each vertex $v$ is in some face, and $v$ is incident to two edges in distinct sets. 

\textbf{Case (b).} $G=V_8$: Using the vertex-numbering in \figref{W}(a), let $E^1:=\{12,34,56,78\}$, $E^2:=\{23,45,67,81\}$, and $E^3:=\{15,26,37,48\}$. Each $E^j$ is a matching of four edges. The claimed properties follow.

\textbf{Case (c).} $G$ is a $(\leq 3)$-sum of two smaller $K_5$-minor-free graphs $G_1$ and $G_2$: Let $C$ be the join set. By induction, there is a partition of each $E(G_i)$ into three sets $E^1_i,E^2_i,E^3_i$ with the desired properties. Permute the set indices so that for each edge $e$ with endpoints in $C$, $e\in E_1^j\cap E_2^j$ for some $j$. This is possible because $C$ is a $(\leq 3)$-clique in $G_1$ and $G_2$. 

For each $j=1,2,3$, let $E^j:=E_1^j\cup E_2^j$. 
If $|C|\leq2$, then  $|E^j|\leq|E_1^j|+|E_2^j|\leq(\V{G_1}-2)+(\V{G_2}-2)=\V{G_1}+\V{G_2}-4\leq\V{G}-2$, as desired. 
Otherwise, $C$ is a triangle in $G_1$ and $G_2$, and $|E_1^j\cap E_2^j|=1$. 
Thus
$|E^j|\leq|E_1^j|+|E_2^j|-1\leq (\V{G_1}-2)+(\V{G_2}-2)-1= \V{G_1}+\V{G_2}-5= \V{G}-2$, as desired. Each triangle of $G$ is in $G_1$ or $G_2$, and thus has one edge in each set $E^j$. 

Consider a $V_8$ subgraph $H$ of $G$. Since $V_8$ is edge-maximal $K_5$-minor-free, $H$ is an induced subgraph. Since $V_8$ is $3$-connected and triangle-free, $H$ is a subgraph of $G_1$ or $G_2$. Thus $H\cap E^j$ is a perfect matching of $H$ by induction.

If $G$ is edge-maximal, then $G_1$ and $G_2$ are both edge-maximal. Thus every vertex $v$ of $G$ is incident to at least two edges in distinct sets (since the same property holds for $v$ in $G_1$ or $G_2$). 
\end{proof}

%%%%%%%%%%%%%%%%%%%%%%%%%%%%%%%%%%%%%%%

For a set $E$ of edges in a graph $G$, a vertex $v$ of $G$ is \emph{$E$-isolated} if $v$ is incident to no edge in $E$.

\begin{lemma}
\lemlabel{OmegaDecompK5}
Suppose that $E$ is a set of edges in a $K_5$-minor-free graph $G$ such that every triangle of $G$ has exactly one edge in $E$, and if $S$ is a subgraph of $G$ isomorphic to $V_8$ then $E(S)\cap E$ is a perfect matching in $S$. Let $V$ be the set of $E$-isolated vertices in $G$. Then $G$ has a planar $\omega$-decomposition $D$ of width $2$ with $V(D)=\{\{v\}:v\in V\}\cup\{\{v,w\}:vw\in E\}$ with no duplicate bags.
\end{lemma}

\begin{proof} 
By \thmref{CharacterisationK5}, we need only consider the following four cases.

\textbf{Case (a).} $G=K_4$: Say $V(G)=\{v,w,x,y\}$. Without loss of generality, $E=\{vw,xy\}$. Thus $V=\emptyset$. Then $D:=K_2$, with bags $\{v,w\}$ and $\{x,y\}$, is the desired decomposition of $G$. Now assume that $G\ne K_4$.

\textbf{Case (b).} $G=V_8$: Thus $E$ is a perfect matching of $G$. Then $D:=K_4$, with one bag for each edge in $E$, is the desired decomposition of $G$.

\textbf{Case (c).} $G$ is planar and has no separating triangle (see \figref{OmegaDecomp}): Fix a plane drawing of $G$. Thus every triangle of $G$ is a face. Initially, let $D$ be the planar decomposition of $G$ with $V(D):=\{\{v\}:v\in V(G)\}$ and $E(D):=\{\{v\}\{w\}:vw\in E(G)\}$. For each edge $vw\in E$, introduce a new bag $\{v,w\}$ in $D$, and replace the edge $\{v\}\{w\}$ by the path  $\{v\}\{v,w\}\{w\}$. Thus $D$ is a planar subdivision of $G$. Now consider each triangle $uvw$ of $G$. Without loss of generality, $vw\in E$. Replace the path $\{v\}\{u\}\{w\}$ in $D$ by the edge $\{u\}\{v,w\}$. Since $uvw$ is a face with only one edge in $E$, $D$ remains planar. Moreover, $D(v)$ is a connected star for each vertex $v$ of $G$. Since $G$ has no separating triangle and $G\ne K_4$, each clique is a $(\leq 3)$-clique. Thus, by construction, each clique is contained in a bag of $D$, or is contained in the union of two adjacent bags of $G$. Therefore $D$ is a planar $\omega$-decomposition of $G$ with width $2$. The order is $\V{G}+|E|$. For each vertex $v\not\in V$, there is an edge incident to $v$ that is in $E$. Choose such an edge $vw\in E$. Thus $\{v\}\{v,w\}$ is an edge of $D$. Contract this edge and merge the bags. By \lemref{ContractDecomp}, $D$ remains a planar $\omega$-decomposition. Now $V(D)=\{\{v\}:v\in V\}\cup\{\{v,w\}:vw\in E\}$.

\Figure{OmegaDecomp}{\includegraphics{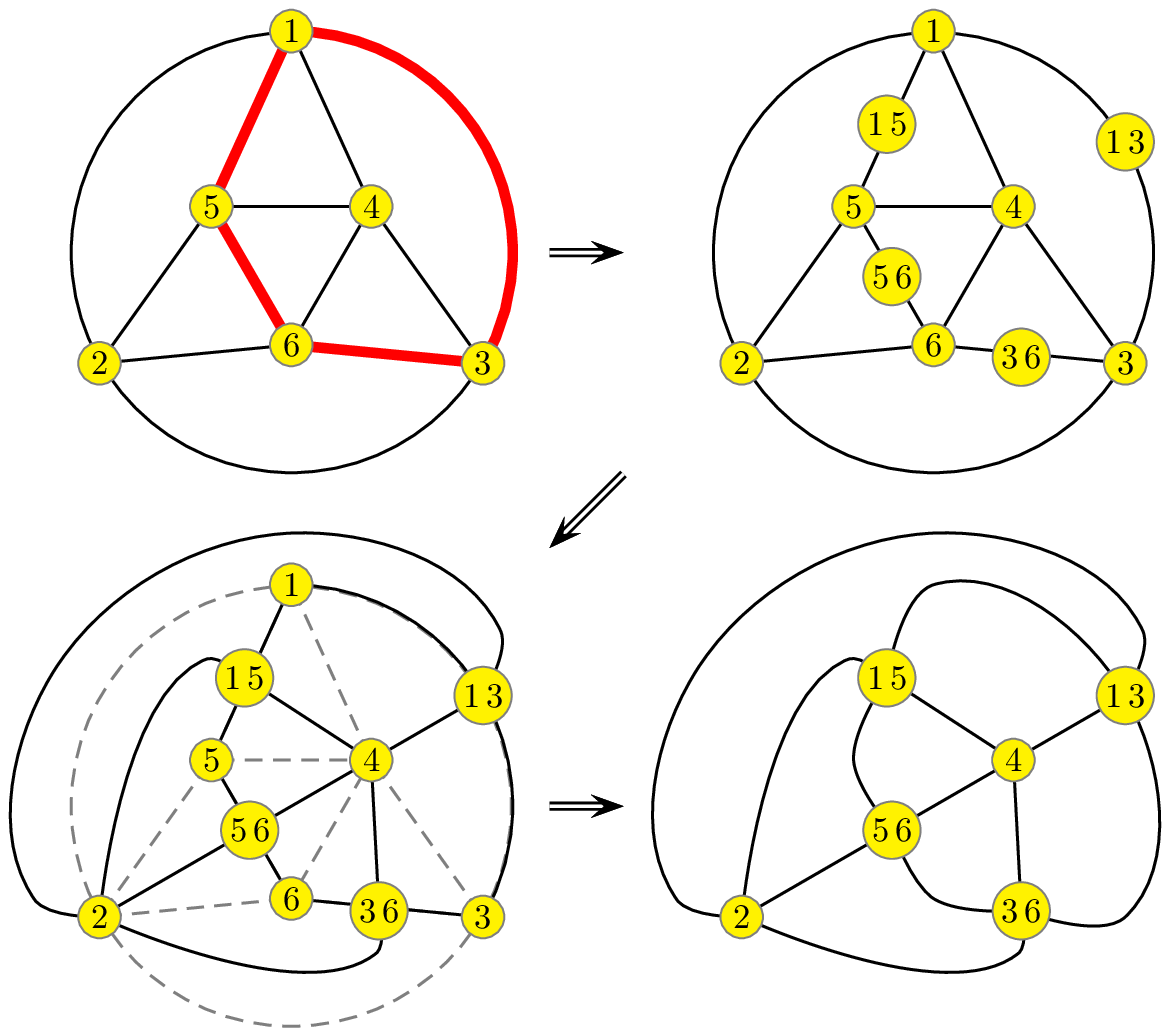}}{An example of the construction of a planar $\omega$-decomposition in case (c) of \lemref{OmegaDecompK5}, where $E=\{13,15,56,36\}$ and $V=\{2,4\}$.}

\textbf{Case (d).} $G$ is a $(\leq 3)$-sum of two smaller $K_5$-minor-free graphs $G_1$ and $G_2$: Let $C$ be the join set. Let $E_1:= E\cap E(G_1)$ and $E_2:= E\cap E(G_2)$. Then every triangle of $G_i$ has exactly one edge in $E_i$. Let $V_i$ be the set of vertices of $G_i$ that are $E_i$-isolated. By induction, each $G_i$ has a planar $\omega$-decomposition $D_i$ of width $2$ with $V(D_i):=E_i\cup V_i$. By \lemref{Sum} with $p=q=4$ and $k_1=k_2=2$, $G$ has a planar $\omega$-decomposition $D$ of width $2$ with $V(D):=V_1\cup E_1\cup V_2\cup E_2$. Moreover, each bag of $D_1$ that intersects $C$ is adjacent to each bag of $D_2$ that intersects $C$. 

If there is a duplicate bag $X$ in $D$, then one copy of $X$ is from $D_1$ and the other copy is from $D_2$, and $X$ intersects $C$. Thus the two copies are adjacent. Contract the edge $XX$ in $D$ into the bag $X$. By  \lemref{ContractDecomp}, $D$ remains a planar $\omega$-decomposition of $G$ with width $2$. Now $D$ has no duplicate bags.   

Every bag in $D$ is either a vertex or an edge of $G$. If $\{v,w\}$ is a bag of $D$, then $vw$ is in $E_1\cup E_2=E$. Conversely, if $vw\in E$ then $vw\in E_1$ or $vw\in E_2$. Thus the bags of cardinality $2$ in $D$ are in one-to-one correspondence with edges in $E$. 

Suppose there is bag $\{v\}$ in $D$ but $v\not\in V$. Then $v\in C$, $v\in V_2$, and $v\not\in V_1$ (or symmetrically, $v\in V_1$ and $v\not\in V_2$). That is, $v$ is incident to an edge $vw\in E_1$ but $v$ is $E_2$-isolated. Now $\{v,w\}$ is a bag in $D_1$, and $\{v\}$ is a bag in $D_2$. These bags are adjacent in $D$. Contract the edge $\{v,w\}\{v\}$ in $D$ into the bag $\{v,w\}$. By \lemref{ContractDecomp}, $D$ remains an $\omega$-decomposition of $G$ with width $2$. An analogous argument applies if $v\in V_1$ and $v\not\in V_2$. We have proved that if $\{v\}$ is a bag of $D$, then $v\in V$. Conversely, if $v\in V$ then $v\in V_1$ or $v\in V_2$ (possibly both), and there is a bag $\{v\}$ in $D$. Thus the singleton bags of $D$ are in one-to-one correspondence with vertices in $V$. 

Therefore $V(D)=\{\{v\}:v\in V\}\cup\{\{v,w\}:vw\in E\}$ with no duplicate bags.
\end{proof}

\begin{proof}[Proof of \thmref{PlanarOmegaDecompK5}]
Let $H$ be an edge-maximal $K_5$-minor-free graph containing $G$ as a spanning subgraph. By \lemref{EdgePartitionK5}, there is a partition of $E(H)$ into three sets $E^1,E^2,E^3$ with the stated properties. We now construct the desired planar $\omega$-decomposition of $H$, which is also the desired decomposition of $G$. Let $V^j$ be the set of $E^j$-isolated vertices in $H$. Consider an $E^1$-isolated vertex $v$. Each edge incident to $v$ is in $E^2\cup E^3$. Since $v$ is incident to (at least) two edges in distinct sets, $v$ is not $E^2$-isolated and is not $E^3$-isolated. In general, $V^i\cap V^j=\emptyset$ for distinct $i$ and $j$\footnote{Note that it is possible for $V_1,V_2,V_3$ to partition $V(H)$; for example, when $H$ is a planar Eulerian triangulation.}. Thus $|V^i|\leq\third\V{H}$ for some $i$. By \lemref{OmegaDecompK5}, $H$ has a planar $\omega$-decomposition of width $2$, with $|E^i|=\V{H}-2$ bags of cardinality $2$, and $|V^i|\leq\frac{1}{3}\V{H}$ bags of cardinality $1$. 
\end{proof}

\begin{proof}[Proof of \thmref{CrossingsK5}]
By \thmref{PlanarOmegaDecompK5}, $G$ has a planar $\omega$-decomposition $D$ of width $2$, with at most $\V{G}-2$ bags of cardinality $2$, and at most
$\frac{1}{3}\V{G}$  bags of cardinality $1$. By \lemref{DecompToCrossing},
\begin{align*}
\CR{G}
\;\leq\;
2\,\Delta(G)^2\!\!\!\sum_{X\in V(D)}\!\!\!\!\!\!\tbinom{|X|+1}{2}
\;\leq\;
2\,\Delta(G)^2\,
\big(3(\V{G}-2)+\tfrac{1}{3}\V{G}\big)
\;<\;
\tfrac{20}{3}\,\Delta(G)^2\,\V{G}
\enspace.
\end{align*}
\end{proof}

%%%%%%%%%%%%%%%%%%%%%%%%%%%%%%%%%%%%%%%%%%%%%%%%%%%%%%%%%%%%%%%%%%%%%%%%%%%

The following two propositions, while not used to prove bounds on the crossing number, are of independent interest. First we consider strong planar $3$-decompositions of $K_5$-minor-free graphs.

\begin{proposition}
\proplabel{Strong3DecompK5}
Every $K_5$-minor-free graph $G$ with $\V{G}\geq3$ has a strong planar $3$-decomposition of width $3$ and order $3\V{G}-8$. Moreover, for all $n\geq3$, there is a planar graph $G$, such that $\V{G}=n$ and every strong planar $3$-decomposition of $G$ with width $3$ has order at least $3\V{G}-8$. 
\end{proposition}

\begin{proof}
Add edges to $G$ so that it is edge-maximal with no $K_5$-minor. This does not affect the claim. By \thmref{CharacterisationK5}, we need only consider the following three cases.

\textbf{Case (a).} $G$ is a planar triangulation with no separating triangle: Let $D$ be the dual graph of $G$. That is, $V(D):=F(G)$ and $E(D):=\{XY:X,Y\in F(G),|X\cap Y|=2\}$. Then $D$ is a planar graph. For each vertex $v$ of $G$, $D(v)$ is the connected cycle consisting of the faces containing $v$. Every $(\leq 3)$-clique of $G$ is contained in some face of $G$ (since $G$ has no separating triangle) and is thus in some bag of $D$. Thus $D$ is a strong planar $3$-decomposition of $G$. The order is $|F(G)|=2\V{G}-4$, which is at most $3\V{G}-8$ unless $\V{G}=3$, in which case one bag suffices.

\textbf{Case (b).} $G=V_8$: Then $\CR{V_8}=1$, as illustrated in \figref{W}(b). Thus by \lemref{CrossingToDecomp}, $G$ has a strong planar decomposition of width $2$ and order $\ceil{\frac{0}{2}}+12+1=13<3\cdot\V{G}-8$; see \figref{W}(d). This decomposition is also a strong planar $3$-decomposition since $K_3\not\subseteq V_8$. 

\textbf{Case (c).} $G$ is a $(\leq 3)$-sum of two smaller $K_5$-minor-free graphs $G_1$ and $G_2$, each with $\V{G_i}\geq3$. Let $C$ be the join set. Thus $C$ is a clique in $G_1$ and in $G_2$. By induction, each $G_i$ has a strong planar $3$-decomposition of width $3$ and order $3\V{G_i}-8$. By \lemref{Sum} with $p=q=k_1=k_2=3$, $G$ has a strong planar $3$-decomposition $D$ of width $3$ and order $3\V{G_1}-8+3\V{G_2}-8=3(\V{G_1}+\V{G_2})-16$. If $|C|\leq 2$ then $\V{G_1}+\V{G_2}\leq \V{G}+2$ and $D$ has order 
at most $3\V{G}-10$. Otherwise $|C|=3$. Since the decompositions of $G_1$ and $G_2$ are strong, $C$ is a bag in both decompositions. Thus $CC$ is an edge of $D$, which can be contracted by \lemref{ContractDecomp}. Therefore $D$ has order $3(\V{G_1}+\V{G_2})-16-1=3(\V{G}+3)-16-1=3\V{G}-8$.

This completes the proof of the upper bound. It remains to prove the lower bound. Observe that in a strong planar $3$-decomposition of width $3$, every triangle is a distinct bag. The lower bound follows since there is a planar graph with $n\geq3$ vertices and $3n-8$ triangles \cite{Wood-PlanarCliques}.
\end{proof}

%%%%%%%%%%%%%%%%%%%%%%%%%%%%%%%%%%%%%%

It follows from Euler's Formula and \thmref{CharacterisationK5} that every $K_5$-minor-free graph $G$ has at most $3\V{G}-6$ edges, and is thus $5$-degenerate (also see \lemref{GenusDegen} below). Thus by \lemref{Degen}, $G$ has a strong $\omega$-decomposition isomorphic to $G$ of width $5$. Since $\omega(G)\leq4$, it is natural to consider strong $\omega$-decompositions of width $4$. 

\begin{proposition}
\proplabel{StrongOmegaDecompK5}
Every $K_5$-minor-free graph $G$ with $\V{G}\geq4$ has a strong planar $\omega$-decomposition of width $4$ and order at most $\frac{4}{3}\V{G}-4$. Moreover, for all $n\geq1$, there is a planar graph $G_n$, such that $\V{G_n}=3n$ and every strong $\omega$-decomposition of $G_n$ with width $4$ has order at least $\frac{7}{6}\V{G_n}-3$. 
\end{proposition}

\begin{proof}
Add edges to $G$ so that it is edge-maximal with no $K_5$-minor. This does not affect the claim. By \thmref{CharacterisationK5}, we need only consider the following three cases.

\textbf{Case (a).} $G$ is a planar triangulation with no separating triangle: If $G=K_4$ then the decomposition with one bag containing all four vertices satisfies the requirements. $K_5$ minus an edge is the only $5$-vertex planar triangulation, and it has a separating triangle. Thus we can assume that $\V{G}\geq6$. By \lemref{EdgePartitionK5}\footnote{The full strength of the Four-Colour Theorem (used in the proof of \lemref{EdgePartitionK5}) is not needed here. That $G$ has (one) set $S$ of $\V{G}-2$ edges such that every face of $G$ has exactly one edge in $S$ quickly follows from Petersen's Matching Theorem applied to the dual; see \cite{BCGMW-SODA06}.}, $G$ has a set $S$ of $\V{G}-2$ edges such that every face of $G$ has exactly one edge in $S$. Each edge $e=vw\in S$ is in two faces $xvw$ and $yvw$; let  $P(e):=\{v,w,x,y\}$. Let $D$ be the graph with $V(D):=\{\{v\}:v\in V(G)\}\cup \{P(e):e\in S\}$, where $\{v\}$ is adjacent to $P(e)$ if and only if $v\in P(e)$. Then $D$ is a planar bipartite graph. For each vertex $v$ of $G$, $D(v)$ is a (connected) star rooted at $\{v\}$. Every face (and thus every triangle) is in some bag of $D$. Thus $D$ is a strong $\omega$-decomposition of $G$ with width $4$ and order $2\V{G}-2$. For each vertex $v$ of $G$, select one bag $P$ containing $v$, and contract the edge $\{v\}P$ in $D$. By \lemref{ContractDecomp}, we obtain a strong $\omega$-decomposition of $G$ with width $4$ and order $\V{G}-2\leq \frac{4}{3}\V{G}-4$ (since $\V{G}\geq 6$).

\textbf{Case (b).} $G=V_8$: As illustrated in \figref{W}(e), $G$ has a strong planar $\omega$-decomposition of width $4$ and order $4<\frac{4}{3}\V{G}-4$. 

\textbf{Case (c).} $G$ is a $(\leq 3)$-sum of two smaller $K_5$-minor-free graphs $G_1$ and $G_2$, each with $\V{G_i}\geq 4$: Let $C$ be the join set. By induction, each $G_i$ has a strong planar $\omega$-decomposition $D_i$ of width $4$ and order at most $\frac{4}{3}\V{G_i}-4$. By \lemref{Sum} with $p\leq 3$ and $q=k_1=k_2=4$, $G$ has a strong planar $\omega$-decomposition of width $4$ and order $\V{D_1}+\V{D_2}\leq \frac{4}{3}\V{G_1}-4+\frac{4}{3}\V{G_2}-4
=\frac{4}{3}(\V{G_1}+\V{G_2})-8 \leq\frac{4}{3}\V{G}-4$.

This completes the proof of the upper bound. It remains to prove the lower bound. Let $G_1:=K_3$. As illustrated in \figref{Example}, construct $G_{n+1}$ from $G_n$ as follows. Insert a triangle inside some face of $G_n$, and triangulate so that each of the three new vertices have degree $4$. This creates seven new triangles and no $K_4$. Thus $\V{G_i}=3n$ and $G_n$ has $7n-6$ triangles. In a strong $\omega$-decomposition of $G_n$, each triangle is in a bag. Since $G_n$ contains no $K_4$, each bag of width $4$ can accommodate at most two triangles. Thus the number of bags is at least half the number of triangles, which is $\frac{7}{2}n-3=\frac{7}{6}\V{G_i}-3$.
\end{proof}

\Figure{Example}{\includegraphics{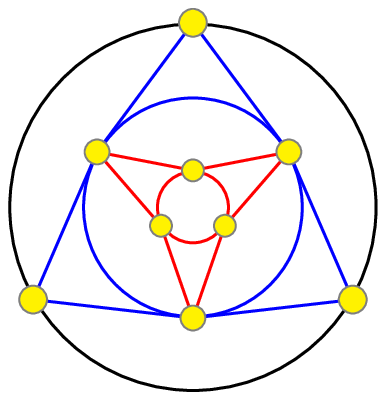}}{The graph $G_3$ in the lower bound of \propref{StrongOmegaDecompK5}.}

We conjecture that for $n\geq2$, the graph $G_n$ in \propref{StrongOmegaDecompK5} actually requires at least $\frac{4}{3}\V{G_n}-4$ bags in every strong $\omega$-decomposition of width $4$.

%%%%%%%%%%%%%%%%%%%%%%%%%%%%%%%%%%%%%%%%%%%%%%%%%%%%%%%%%%%%%%%%%%%%%%%%%%%%%%%%
\section{Graphs Embedded on a Surface}
\seclabel{Surfaces}
%%%%%%%%%%%%%%%%%%%%%%%%%%%%%%%%%%%%%%%%%%%%%%%%%%%%%%%%%%%%%%%%%%%%%%%%%%%%%%%%

Recall that \SSS\ is the orientable surface with $\gamma$ handles. As illustrated in \figref{Cycles}, a \emph{cycle} in \SSS\ is a closed curve in the surface. A cycle is \emph{contractible} if it is contractible to a point in the surface. A noncontractible cycle is \emph{separating} if it separates \SSS\ into two connected components. 

\Figure{Cycles}{\includegraphics{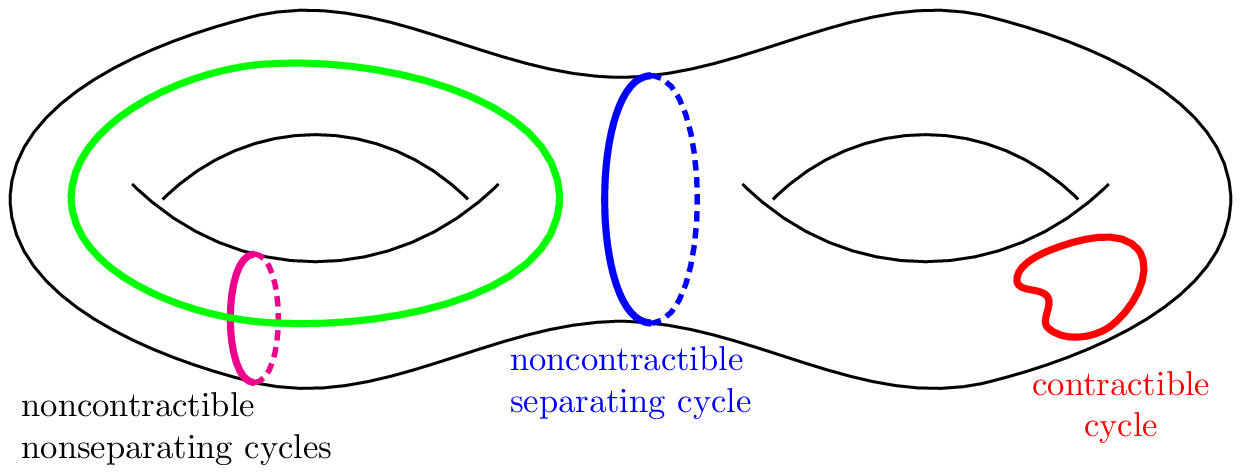}}{Cycles on the double torus.}

Let $G$ be a graph embedded in \SSS. A \emph{noose} of $G$ is a cycle $C$ in \SSS\ that does not intersect the interior of an edge of $G$. Let $V(C)$ be the set of vertices of $G$ intersected by $C$. The \emph{length} of $C$ is $|V(C)|$. 
Pach and T\'{o}th~\cite{PachToth-GD05} proved that, for some constant $c_\gamma$, the crossing number of every graph $G$ of genus $\gamma$ satisfies\footnote{Recently Djidjev,and Vr{\v{t}}o~\cite{DjidVrto-ICALP06} greatly improved the constant $c_\gamma$ in \eqnref{PachToth}, and B{\"o}r{\"o}czky~et~al.~\cite{BPT-IJFCS06} extended the bound to graphs embedded on non-orientable surfaces.}
\begin{equation}
\eqnlabel{PachToth}
\CR{G}\;\leq\; c_\gamma\!\!\!\sum_{v\in V(G)}\!\!\!\deg(v)^2
\;\leq\; 2c_\gamma\,\Delta(G)\,\E{G}\enspace.
\end{equation}

The following lemma is probably well known.  

\begin{lemma}
\lemlabel{NumberEdgesGenus}
Every graph with $n\geq 3$ vertices and genus $\gamma$ has at most $(\sqrt{3\gamma}+3)n-6$ edges.
\end{lemma}

\begin{proof}
Let $G$ be an edge-maximal graph with $n$ vertices, $m$ edges, and genus $\gamma$. Suppose that an embedding of $G$ in \SSS\ has $f$ faces. Euler's Formula states that $m-n-f=2\gamma-2$. Since $G$ is edge-maximal, every face is bounded by three edges and each edge is in the boundary of two faces. Thus $m=3f/2$ and $f=2m/3$. Hence $m=3n+6\gamma-6$. If $\gamma=0$ then we are done. Now assume that $\gamma\geq1$. We need to prove that $3n+6\gamma\leq (\sqrt{3\gamma}+3)n$. That is, $6\gamma\leq\sqrt{3\gamma}n$, or equivalently $\gamma\leq n^2/12$, which is true since $K_n$ has genus $\ceil{(n-3)(n-4)/12}<n^2/12$; see \cite{Ringel74, MoharThom}.
\end{proof}

\Eqnref{PachToth} and \lemref{NumberEdgesGenus} imply that 
\begin{equation}
\eqnlabel{PachTothPlus}
\CR{G}\;\leq\; c_\gamma\,\Delta(G)\,\V{G}\enspace.
\end{equation}

By \lemref{CrossingToDecomp}, $G$ has a planar decomposition of width $2$ and order $c_\gamma\,\Delta(G)\,\V{G}$. We now provide an analogous result without the dependence on $\Delta(G)$, but at the expense of an increased bound on the width.

\begin{theorem}
\thmlabel{GenusPlanarDecomp}
Every graph $G$ with genus $\gamma$ has a planar decomposition of width $2^{\gamma}$ and order $3^{\gamma}\V{G}$. 
\end{theorem}

The key to the proof of \thmref{GenusPlanarDecomp} is the following lemma, whose proof is inspired by similar ideas of Pach and T\'{o}th~\cite{PachToth-GD05}.

\begin{lemma}
\lemlabel{ReduceGenus}
Let $G$ be a graph with a $2$-cell embedding in \SSS\ for some $\gamma\geq1$. Then $G$ has a decomposition of width $2$, genus at most $\gamma-1$, and order $3\V{G}$.
\end{lemma}

\begin{proof}
Since $\gamma\geq1$, \SSS\ has a noncontractible nonseparating cycle, which can be continuously deformed into a noncontractible nonseparating noose in $G$. Let $C$ be a noncontractible nonseparating noose in $G$ of minimum length $k:=|V(C)|$. Orient $C$ and let $V(C):=(v_1,v_2,\dots,v_k)$ in the order around $C$. For each vertex $v_i\in V(C)$, let $E^\ell(v_i)$ and $E^r(v_i)$ respectively be the set of edges incident to $v_i$ that are on the left-hand side and right-hand side of $C$ (with respect to the orientation). Cut the surface along $C$, and attach a disk to each side of the cut. Replace each vertex $v_i\in V(C)$ by two vertices $v_i^\ell$ and $v_i^r$ respectively incident to the edges in $E^\ell(v_i)$ and $E^r(v_i)$. Embed $v_i^\ell$ on the left-hand side of the cut, and embed $v^r$ on the right-hand side of the cut. We obtain a graph $G'$ embedded in a surface of genus at most $\gamma-1$ (since $C$ is nonseparating). 

Let $L:=\{v_i^\ell:v\in V(C)\}$ and $R:=\{v_i^r:v\in V(C)\}$. By Menger's Theorem, the maximum number of disjoint paths between $L$ and $R$ in $G'$ equals the minimum number of vertices that separate $L$ from $R$ in $G'$. Let $Q$ be a minimum set of vertices that separate $L$ from $R$ in $G'$. Then there is a noncontractible nonseparating noose in $G$ that only intersects vertices in $Q$. (It is nonseparating in $G$ since $L$ and $R$ are identified in $G$.)\ Thus $|Q|\geq k$ by the minimality of $|V(C)|$. Hence there exist $k$ disjoint paths $P_1,P_2,\dots,P_k$ between $L$ and $R$ in $G'$, where the endpoints of $P_i$ are $v_i^\ell$ and $v_{\sigma(i)}^r$, for some permutation $\sigma$ of $[1,k]$. In the disc with $R$ on its boundary, draw an edge from each vertex $v_{\sigma(i)}^r$ to $v_i^r$ such that no three edges cross at a single point and every pair of edge cross at most once.  Add a new vertex $x_{i,j}$ on each crossing point between edges $v_{\sigma(i)}^rv_i^r$ and $v_{\sigma(j)}^rv_j^r$.  Let $G''$ be the graph obtained. Then $G''$ is embedded in $S_{\gamma-1}$. 

We now make $G''$ a decomposition of $G$. Replace $v_i^\ell$ by $\{v_i\}$ and replace $v_i^r$ by $\{v_i\}$. Replace every other vertex $v$ of $G$ by $\{v\}$.
Replace each `crossing' vertex $x_{i,j}$ by $\{v_i,v_j\}$. Now for each vertex $v_i\in V(C)$, add $v_i$ to each bag on the path $P_i$ from $v_i^\ell$ to $v_{\sigma(i)}^r$. Thus $G''(v_i)$ is a (connected) path. Clearly $G''(v)$ and $G''(w)$ touch for each edge $vw$ of $G$. Hence $G''$ is a decomposition of $G$ with genus at most $\gamma-1$. Since the paths $P_1,P_2,\dots,P_k$ are pairwise disjoint, the width of the decomposition is $2$. 

It remains to bound the order of $G''$. Let $n:=\V{G}$. Observe that $G''$ has at most $n+k+\binom{k}{2}$ vertices. One of the paths $P_i$ has at most $\frac{n+k}{k}$ vertices. For ease of counting, add a cycle to $G'$ around $R$. Consider the path in $G'$ that starts at $v_i^\ell$, passes through each vertex in $P_i$, and then takes the shortest route from $v_{\sigma(i)}^r$ around $R$ back to $v_i^r$. The distance between $v_{\sigma(i)}^r$ and $v_i^r$ around $R$ is at most $\frac{k}{2}$. This path in $G'$ forms a noncontractible nonseparating noose in $G$ (since if two cycles in a surface cross in exactly one point, then both are noncontractible).

The length of this noose in $G$ is at most $\frac{n+k}{k}-1+\frac{k}{2}$ (since $v_i^\ell$ and $v_i^r$ both appeared in the path).  Hence $\frac{n+k}{k}-1+\frac{k}{2}\geq k$ by the minimality of $|V(C)|$. Thus $k\leq\sqrt{2n}$. Therefore $G''$ has at most $n+\sqrt{2n}+\binom{\sqrt{2n}}{2}\leq 3n$ vertices.
\end{proof}

\begin{proof}[Proof of \thmref{GenusPlanarDecomp}]
We proceed by induction on $\gamma$. If $\gamma=0$ then $G$ is planar, and $G$ itself is a planar decomposition of width $1=2^0$ and order $n=3^0n$. Otherwise, by \lemref{ReduceGenus}, $G$ has a decomposition $D$ of width $2$, genus $\gamma-1$, and order $3n$. By induction, $D$ has a planar decomposition of width $2^{\gamma-1}$ and order $3^{\gamma-1}(3n)=3^{\gamma}n$. By \lemref{Composition} with $p=k=2$, and $\ell=2^{\gamma-1}$, $G$ has a planar decomposition of width $2\cdot2^{\gamma-1}=2^{\gamma}$ and order $3^{\gamma}n$.
\end{proof}

\thmref{GenusPlanarDecomp} and \lemref{DecompToCrossing} imply that every graph $G$ with genus $\gamma$ has crossing number  $\CR{G}\leq 12^{\gamma}\,\Delta(G)^2\,\V{G}$, which for fixed $\gamma$, is weaker than the bound of Pach and T\'{o}th~\cite{PachToth-GD05} in \eqnref{PachTothPlus}. The advantage of our approach is that it generalises for graphs with an arbitrary excluded minor (and the dependence on $\gamma$ is much smaller). 

%For applications to follow we need a planar $\omega$-decomposition of a graph with bounded genus.

We now prove that a graph $G$ embedded on a surface has an $\omega$-decomposition with small width and linear order. To do so, we apply \lemref{DegenDegen}, which requires a bound on the degeneracy of $G$. 

\begin{lemma}
\lemlabel{GenusDegen}
Every graph $G$ of genus $\gamma$ is $(2\sqrt{3\gamma}+6)$-degenerate.
If $\sqrt{3\gamma}$ is an integer then $G$ is $(2\sqrt{3\gamma}+5)$-degenerate.
\end{lemma}

\begin{proof}
By \lemref{NumberEdgesGenus}, $G$ has average degree $\frac{2\E{G}}{\V{G}}<2(\sqrt{3\gamma}+3)$. Thus $G$ has a vertex of degree less than $2\sqrt{3\gamma}+6$. Moreover, if $\sqrt{3\gamma}$ is an integer, then $G$ has a vertex of degree at most $2\sqrt{3\gamma}+5$. The result follows since every subgraph of $G$ has genus at most $\gamma$. 
\end{proof}

\begin{theorem}
\thmlabel{GenusPlanarOmegaDecomp}
Every graph $G$ of genus $\gamma$ has a planar $\omega$-decomposition of width $2^{\gamma}(2\sqrt{3\gamma}+7)$ and order $3^{\gamma}\,\V{G}$.
\end{theorem}

\begin{proof}
By \thmref{GenusPlanarDecomp}, $G$ has a planar decomposition $D$ of width at most $2^{\gamma}$ and order $3^{\gamma}n$. By \lemref{GenusDegen}, $G$ is $(2\sqrt{3\gamma}+6)$-degenerate. Thus by \lemref{DegenDegen}, $G$ has a planar $\omega$-decomposition isomorphic to $D$ with width $2^{\gamma}(2\sqrt{3\gamma}+7)$.
\end{proof}

%%%%%%%%%%%%%%%%%%%%%%%%%%%%%%%%%%%%%%%%%%%%%%%%%%%%%%%%%%%%%%%%%%%%%%%%%%%%%%%%
\section{$H$-Minor-Free Graphs}
\seclabel{Minors}
%%%%%%%%%%%%%%%%%%%%%%%%%%%%%%%%%%%%%%%%%%%%%%%%%%%%%%%%%%%%%%%%%%%%%%%%%%%%%%%%

For integers $h\geq1$ and $\gamma\geq0$, Robertson and Seymour~\cite{RS-GraphMinorsXVI-JCTB03} defined a graph $G$ to be \emph{$h$-almost embeddable} in \SSS\ if $G$ has a set $X$ of at most $h$ vertices such that $G\setminus X$ can be written as $G_0\cup G_1\cup \dots \cup G_h$ such that:
\begin{itemize}
\item $G_0$ has an embedding in $\SSS$,
\item the graphs $G_1,G_2,\dots,G_h$  (called \emph{vortices}) are pairwise disjoint,
\item there are faces\footnote{Recall that we equate a face with the set of vertices on its boundary.} $F_1,F_2,\dots,F_h$ of the embedding of $G_0$ in $\SSS$, such that each $F_i=V(G_0)\cap V(G_i)$,% and $t_i:=|F_i|$,
\item if $F_i=(u_{i,1},u_{i,2},\dots,u_{i,|F_i|})$ in clockwise order about the face, then $G_i$ has a strong $|F_i|$-cycle decomposition $Q_i$ of width $h$, such that each vertex $u_{i,j}$ is in the $j$-th bag of $Q_i$.
\end{itemize}

The following `characterisation' of $H$-minor-free graphs is a deep theorem by Robertson and Seymour~\cite{RS-GraphMinorsXVI-JCTB03}; see the recent survey by Kawarabayashi and Mohar~\cite{KM-GC07}.

\begin{theorem}[Graph Minor Structure Theorem \cite{RS-GraphMinorsXVI-JCTB03}]
\thmlabel{RS}
For every graph $H$ there is a positive integer $h=h(H)$, such that every $H$-minor-free graph $G$ can be obtained by $(\leq h)$-sums of graphs that are $h$-almost embeddable in some surface in which $H$ cannot be embedded.
\end{theorem}

The following theorem is one of the main contributions of this paper. 

\begin{theorem}
\thmlabel{DecompMinorFree}
For every graph $H$ there is an integer $k=k(H)$, such that every $H$-minor-free graph $G$ has a planar $\omega$-decomposition of width $k$ and order $\V{G}$.
\end{theorem}

We prove \thmref{DecompMinorFree} by a series of lemmas.

\begin{lemma}
\lemlabel{AlmostEmbedPlanarOmegaDecomp}
Every graph $G$ that is $h$-almost embeddable in \SSS\ has a planar decomposition of width $h(2^{\gamma}+1)$ and order $3^{\gamma}\,\V{G}$.
\end{lemma}

\begin{proof}
By \thmref{GenusPlanarDecomp}, $G_0$ has a planar decomposition $D$ of width at most $2^{\gamma}$ and order $3^{\gamma}\,\V{G_0}\leq 3^{\gamma}\,\V{G}$. We can assume that $D$ is connected. For each vortex $G_i$, add each vertex in the $j$-th bag of $Q_i$ to each bag of $D$ that contains $u_{i,j}$. The bags of $D$ now contain at most $2^{\gamma}h$ vertices. Now add $X$ to every bag. The bags of $D$ now contain at most $(2^{\gamma}+1)h$ vertices. For each vertex $v$ that is not in a vortex, $D(v)$ is unchanged by the addition of the vortices, and is thus connected. For each vertex $v$ in a vortex $G_i$, $D(v)$ is the subgraph of $D$ induced by the bags (in the decomposition of $G_0$) that contain $u_{i,j}$, where $v$ is in the $j$-th bag of $Q_i$. Now $Q_i(v)$ is a connected subgraph of the cycle $Q_i$, and for each vertex $u_{i,j}$, the subgraphs $G_0(u_{i,j})$ and
$G_0(u_{i,j+1})$ touch. Thus $D(v)$ is connected. (This argument is similar to that used in \lemref{Composition}.)\ $D(v)$ is connected for each vertex $v\in X$ since $D$ itself is connected. 
\end{proof}

\begin{lemma}
\lemlabel{AlmostEmbedDegen}
For all integers $h\geq1$ and $\gamma\geq0$ there is a constant $d=d(h,\gamma)$, such that every graph $G$ that is $h$-almost embeddable in \SSS\ is $d$-degenerate.
\end{lemma}

\begin{proof}
If $G$ is $h$-almost embeddable in \SSS\ then every subgraph of $G$ is
$h$-almost embeddable in \SSS. Thus it suffices to prove that if $G$ has $n$ vertices and $m$ edges, then its average degree $\frac{2m}{n}\leq d$. 
Say each $G_i$ has $m_i$ edges. $G$ has at most $hn$ edges incident to $X$. Thus $m\leq hn + \sum_{i=0}^h m_i$. By \lemref{NumberEdgesGenus}, $m_0<(\sqrt{3\gamma}+3)n$. Now $m_i\leq \binom{h}{2}|F_i|$ by \Eqnref{NumberEdgesStrongDecomp} with $D$ an $|F_i|$-cycle. Since $G_1,G_2,\dots,G_h$ are pairwise disjoint, $\sum_{i=1}^h m_i\leq \binom{h}{2}n$. Thus $m< (h + \sqrt{3\gamma}+3 + \binom{h}{2})n$. Taking $d=h(h+1)+2\sqrt{3\gamma}+6$ we are done. 
\end{proof}

%%%%%%%%%%%%%%%%

\threelemref{DegenDegen}{AlmostEmbedPlanarOmegaDecomp}{AlmostEmbedDegen} imply:

\begin{corollary}
\corlabel{DegenAlmostEmbedPlanarOmegaDecomp}
For all integers $h\geq1$ and $\gamma\geq0$ there is a constant $k=k(h,\gamma)\geq h$, such that every graph $G$ that is $h$-almost embeddable in \SSS\ has a planar $\omega$-decomposition of width $k$ and order $3^{\gamma}\,\V{G}$. \qed
\end{corollary}

Now we bring in $(\leq h)$-sums.

\begin{lemma}
\lemlabel{DecompAlmostEmbedSums}
For all integers $h\geq1$ and $\gamma\geq0$, every graph $G$ that can be obtained by $(\leq h)$-sums of graphs that are $h$-almost embeddable in \SSS\ has a planar $\omega$-decomposition of width $k$ and order $\max\{1,3^{\gamma}(h+1)(\V{G}-h)\}$, where $k=k(h,\gamma)$ from \corref{DegenAlmostEmbedPlanarOmegaDecomp}.
\end{lemma}

\begin{proof}
We proceed by induction on $\V{G}$. If $\V{G}\leq h$ then the decomposition of $G$ with all its vertices in a single bag satisfies the claim (since $k\geq h$).

Now assume that $\V{G}\geq h+1$. If $G$ is $h$-almost embeddable in $\SSS$, then by \corref{DegenAlmostEmbedPlanarOmegaDecomp}, $G$ has a planar $\omega$-decomposition of width $k$ and order $3^{\gamma}\V{G}$, which, since $\V{G}\geq h+1$, is at most $3^{\gamma}(h+1)(\V{G}-h)$, as desired. 

Otherwise, $G$ is a $(\leq h)$-sum of graphs $G_1$ and $G_2$, each of which, by induction, has a planar $\omega$-decomposition of width $k$ and order $\max\{1,3^{\gamma}(h+1)(\V{G_i}-h)\}$. By \lemref{Sum}, $G$ has a planar $\omega$-decomposition $D$ of width $k$ and order 
\begin{equation*}
\V{D}=\max\{1,3^{\gamma}(h+1)(\V{G_1}-h)\}+\max\{1,3^{\gamma}(h+1)(\V{G_2}-h)\}\enspace.
\end{equation*}

Without loss of generality, $\V{G_1}\leq\V{G_2}$. If $\V{G_2}\leq h$ then $\V{D}=2\leq 3^{\gamma}(h+1)(\V{G}-h)$, as desired. If $\V{G_1}\leq h$ and $\V{G_2}\geq h+1$, then $\V{D}=1+3^{\gamma}(h+1)(\V{G_2}-h)$, which, 
since $\V{G}\geq\V{G_2}+1$, is at most $3^{\gamma}(h+1)(\V{G}-h)$, as desired. Otherwise, both $\V{G_1}\geq h+1$ and $\V{G_2}\geq h+1$. Thus the order of $D$ is $3^{\gamma}(h+1)(\V{G_1}+\V{G_2}-2h)\leq 3^{\gamma}(h+1)(\V{G}-h)$, as desired. 
\end{proof}

%%%%%%%%%%%%%%%%%%%%%%%%

\begin{proof}[Proof of \thmref{DecompMinorFree}]
Let $h=h(H)$ from \thmref{RS}. Let \SSS\ be the surface in \thmref{RS} in which $H$ cannot be embedded. If $G$ has no $H$-minor then, by \thmref{RS}, $G$ can be obtained by $(\leq h)$-sums of graphs that are $h$-almost embeddable in \SSS. 
By \lemref{DecompAlmostEmbedSums}, $G$ has a planar $\omega$-decomposition of width $k$ and order $3^{\gamma}(h+1)\V{G}$, where 
$k=k(h,\gamma)$ from \corref{DegenAlmostEmbedPlanarOmegaDecomp}. By \lemref{ReduceOrder}, $G$ has a planar $\omega$-decomposition of width $k'$ and order $\V{G}$, for some $k'$ only depending on $k$, $\gamma$ and $h$ (all of which only depend on $H$).
\end{proof}

\thmref{DecompMinorFree} and \lemref{DecompToCrossing} imply the following quantitative version of \thmref{CrossingMinorFree}.

\begin{corollary}
\corlabel{CrossingMinorFreeQuantative}
For every graph $H$ there is a constant $c=c(H)$, such that every $H$-minor-free graph $G$ has crossing number at most $c\,\Delta(G)^2\,\V{G}$.\qed
\end{corollary}

It is an open problem whether the dependence on $\Delta(G)$ in \corref{CrossingMinorFreeQuantative} can be reduced from quadratic to linear (even with $H=K_5$). We conjecture the stronger result that every $H$-minor-free graph $G$ has crossing number 
\begin{equation*}
\CR{G} \;\leq\; c_H \sum_{v\in V(G)}\!\!\!\deg(v)^2\enspace.
\end{equation*}
Pach and T\'{o}th~\cite{PachToth-GD05} proved this conjecture for graphs of bounded genus.

%%%%%%%%%%%%%%%%%%%%%%%%%%%%%%%%%%%%%%%%%%%%%%%%%%%%%%%%%%%%%
\section{Graph Partitions}
\seclabel{Partitions}
%%%%%%%%%%%%%%%%%%%%%%%%%%%%%%%%%%%%%%%%%%%%%%%%%%%%%%%%%%%%%

A \emph{partition} of a graph is a proper partition of its vertex set. Each part of the partition is called a \emph{bag}. The \emph{width} of partition is the maximum number of vertices in a bag. The \emph{pattern} (or \emph{quotient graph}) of a partition is the graph obtained by identifying the vertices in each bag, deleting loops, and replacing parallel edges by a single edge. Observe that a graph $G$ has a decomposition $D$ of spread $1$ if and only if $G$ has a partition whose pattern is a subgraph of $D$. 

A \emph{tree-partition} is a partition whose pattern is a forest. The \emph{tree-partition-width}\footnote{Tree-partition-width has also been called \emph{strong tree-width} \cite{Seese85,BodEng-JAlg97}.} of a graph $G$ is the minimum width of a tree-partition of $G$, and is denoted by $\tpw{G}$. Tree-partitions were independently introduced by Seese~\cite{Seese85} and Halin~\cite{Halin91}, and have since been investigated by a number of authors \cite{BodEng-JAlg97, Bodlaender-DMTCS99, DO-JGT95, DO-DM96, Edenbrandt86, DMW-SJC05, Wood-JGT06}. 

%OTHER REFERENCES \cite{GeorgeLiu80,Waller76,ER-TCS90,ER-TCS90b}

A graph with bounded degree has bounded tree-partition-width if and only if it has bounded tree-width \cite{DO-DM96}. In particular, Seese~\cite{Seese85} proved the lower bound, 
\begin{equation*}
\label{WidthLowerBound}
2\cdot\tpw{G}\geq\tw{G}+1\enspace,
\end{equation*}
which is tight for even complete graphs. The best known upper bound is
\begin{equation}
\eqnlabel{WidthUpperBound}
\tpw{G}\leq \tfrac{5}{2}\big(\tw{G}+1\big)\big(\tfrac{7}{2}\,\Delta(G)-1\big)\enspace,
\end{equation}
which was obtained by the first author \cite{Wood-TreePartitions} using a minor improvement to a similar result by an anonymous referee of the paper by Ding and Oporowski~\cite{DO-JGT95}. See \cite{ADOV-JCTB03, DHK-FOCS05, DDOSRSV-JCTB04, DOSV-Comb98, DOSV-JCTB00, NesOdM-GradI} for other results related to tree-width and graph partitions.

Here we consider more general types of partitions. A partition is \emph{planar} if its pattern is planar. A relationship between planar partitions and rectilinear drawings is described in the following lemma\footnote{Note that \lemref{ProduceDrawing} bounds the number of crossings per edge; see \cite{PachToth-Comb97} for related results.}. 

\begin{lemma}
\lemlabel{ProduceDrawing}
Every graph $G$ with a planar partition of width $p$ has a rectilinear drawing in which each edge crosses at most $2\,\Delta(G)\,(p-1)$ other edges. Hence
\begin{equation*}
\RCR{G}\leq(p-1)\,\Delta(G)\,\E{G}.
\end{equation*}
\end{lemma}

\begin{proof} 
%Let $H$ be the pattern of a planar partition of $G$ with width $p$. Let $B_x$ be the set of vertices of $G$ in the bag that corresponds to each vertex $x$ of $H$. 
%
%By the F{\'a}ry-Wagner  Theorem, $H$ has a rectilinear drawing with no crossings. Let $\epsilon>0$. Let $D_\epsilon(x)$ be the disc of radius $\epsilon$ centred at each vertex $x$ of $H$. For each edge $xy$ of $H$, let $D_\epsilon(xy)$ be the union of all line-segments with one endpoint in $D_\epsilon(x)$ and one endpoint in $D_\epsilon(y)$. For some  $\epsilon>0$, we have $D_\epsilon(x)\cap D_\epsilon(y)=\emptyset$ for all distinct vertices $x$ and $y$ of $H$, and $D_\epsilon(xy)\cap D_\epsilon(pq)=\emptyset$ for all edges $xy$ and $pq$ of $H$ that have no endpoint in common. 
%
%Position each vertex $v\in B_x$ inside $D_\epsilon(x)$ so that the vertices of $G$ are in general position (no three collinear). Draw every edge of $G$ straight. Thus no edge passes through a vertex. 
%
Apply the construction from \lemref{DecompToCrossing} with $s(v)=1$ for every vertex $v$. We obtain a rectilinear drawing of $G$. Consider an edge $vw$ of $G$. Say $v$ is in bag $X$, and $w$ is in bag $Y$. Then $vw$ is drawn inside $D_\epsilon(XY)$. Thus, if two edges $e_1$ and $e_2$ of $G$ cross, then an endpoint of $e_1$ and an endpoint of $e_2$ are in a common bag, and $e_1$ and $e_2$ have no endpoint in common. Thus each edge of $G$ crosses at most $2\,\Delta(G)\,(p-1)$ other edges, and $\RCR{G}\leq\half\sum_e2\,\Delta(G)\,(p-1)=\Delta(G)\,(p-1)\E{G}$.
\end{proof}

A graph is \emph{outerplanar} if it has a plane drawing with all the vertices on the outerface. Obviously, $\CCR{G}=0$ if and only if $G$ is outerplanar.  A partition is \emph{outerplanar} if its pattern is outerplanar. 

\begin{lemma}
\lemlabel{ProduceConvexDrawing}
Every graph $G$ with an outerplanar partition of width $p$ has a convex drawing in which each edge crosses at most $2\,\Delta(G)\,(p-1)$ other edges. Hence \begin{equation*}
\CCR{G}\leq(p-1)\,\Delta(G)\,\E{G}\enspace.
\end{equation*}
\end{lemma}

\begin{proof}
Apply the proof of \lemref{ProduceDrawing} starting from a plane convex drawing of the pattern.
\end{proof}

Since every forest is outerplanar, every graph $G$ has an outerplanar partition of width \tpw{G}. Thus \lemref{ProduceConvexDrawing} and \Eqnref{WidthUpperBound} imply the following quantitative version of \thmref{ConvexCrossingTreewidth}.

\begin{corollary}
\corlabel{TreewidthCrossings}
Every graph $G$ has a convex drawing in which each edge crosses less than
\begin{equation*}
5\,\Delta(G)\,\big(\tw{G}+1\big)\big(7\,\Delta(G)-1\big)
\end{equation*}
other edges. Hence 
\begin{equation*}
\CCR{G}
<\tfrac{17}{2}\big(\tw{G}+1\big)\,\Delta(G)^2\,\E{G}
<\tfrac{17}{2}\,\tw{G}\,\big(\tw{G}+1\big)\,\Delta(G)^2\,\V{G}
\enspace.
\end{equation*}\qed
\end{corollary}

Alon~et~al.~\cite{AST-JAMS90} proved that every $H$-minor free graph $G$ has tree-width at most $c(H)\,\sqrt{\V{G}}$ for some constant $c(H)$; also see \cite{Grohe-Comb03, DemHaj-SODA05}. Thus \corref{TreewidthCrossings} implies: 

\begin{corollary}
For every graph $H$ there is a constant $c=c(H)$, such that every graph $G$ with no $H$-minor has a convex drawing in which each edge crosses less than
$c\,\Delta(G)^2\,\sqrt{\V{G}}$ other edges. Hence 
\begin{equation*}
\CCR{G}
  <c\,\Delta(G)^2\,\sqrt{\V{G}}\,\E{G}
\;<c\;\Delta(G)^2\,\V{G}^{3/2}
\enspace.
\end{equation*}\qed
\end{corollary}

%%%%%%%%%%%%%%%%%%%%%%%%%%%%%%

Note the following result which is converse to \corref{TreewidthCrossings}.

\begin{proposition}
\proplabel{ConvexTreewidth}
Suppose that a graph $G$ has a convex drawing such that whenever two edges $e$ and $f$ cross, $e$ or $f$ crosses at most $k$ edges. Then $G$ has tree-width $\tw{G}\leq 3k+11$.
\end{proposition}

\begin{proof}
First we construct a strong planar decomposition $D$ of $G$ (in a similar way to the proof of \lemref{CrossingToDecomp}). Replace each vertex $v$ of $G$ by the bag $\{v\}$ in $D$. Orient each edge of $G$. Replace each crossing between arcs $(v,w)$ and $(x,y)$ of $G$ by the bag $\{v,x\}$ in $D$. For each arc $(v,w)$ of $G$,  for some vertex $x$ of $G$, there is an edge $\{v,x\}\{w\}$ in $D$; replace this edge by the path $\{v,x\}\{v,w\}\{w\}$. Thus $D$ is a strong planar decomposition of $G$ with width $2$. Observe that the distance between each bag in $D$ and some bag $\{v\}$ on the outerface is at most $\floor{\frac{k}{2}}+1$. Thus $D$ is $(\floor{\frac{k}{2}}+2)$-outerplanar\footnote{An outerplanar graph is called \emph{$1$-outerplanar}. A plane graph is \emph{$k$-outerplanar} if the graph obtained by deleting the vertices on the outerface is ($k-1)$-outerplanar.}. Bodlaender~\cite{Bodlaender88} proved that every $d$-outerplanar graph has tree-width at most $3d-1$. Thus $D$ has tree-width at most $3\floor{\frac{k}{2}}+5$. That is, some tree $T$ is a strong decomposition of $D$ with width at most $3\floor{\frac{k}{2}}+6$. By \lemref{Composition} with $p=2$, $G$ has a strong decomposition isomorphic to $T$ with width at most $6\floor{\frac{k}{2}}+12$. That is, $G$ has tree-width at most $6\floor{\frac{k}{2}}+11$.
\end{proof}

%%%%%%%%%%%%%%%%%%%%%%%%%%%%%%%%%%%%%%%%%%%%%%%%%%%%%%%%%%%%%%%%%%%%%%%%%%%%%%%%
\section{$K_{3,3}$-Minor-Free Graphs}
\seclabel{K33}
%%%%%%%%%%%%%%%%%%%%%%%%%%%%%%%%%%%%%%%%%%%%%%%%%%%%%%%%%%%%%%%%%%%%%%%%%%%%%%%%

In this section we prove \thmref{RectilinearCrossingK33}, which gives an upper bound on the rectilinear crossing number of $K_{3,3}$-minor-free graphs. The proof employs the following characterisation by Wagner~\cite{Wagner37}. 

\begin{lemma}[\cite{Wagner37}]
\lemlabel{CharacterisationK33}
A graph $G$ is $K_{3,3}$-minor-free if and only if $G$ can be obtained by $(\leq 2)$-sums from planar graphs and $K_5$.
\end{lemma}

\begin{lemma}
Let $G$ be a $K_{3,3}$-minor-free graph. For every edge $e$ of $G$, there is a matching $M$ in $G$ with the following properties: 
\begin{itemize}
\item $|M|\leq \third(\V{G}-2)$,
\item each edge in $M$ is disjoint from $e$,
\item contracting $M$ gives a planar graph.
\end{itemize}
\end{lemma}

\begin{proof}
If $G$ is planar, then the lemma is satisfied with $M=\emptyset$. Suppose that $G=K_5$. Let $vw$ be an edge of $G$ that is disjoint from $e$. Let $M:=\{vw\}$. Then $|M|=1=\third(\V{G}-2)$. The graph obtained by contracting $vw$ is $K_4$, which is planar. 

Now assume that $G$ is not planar and not $K_5$. By \lemref{CharacterisationK33}, $G$ is a $(\leq 2)$-sum of two smaller $K_{3,3}$-minor-free graphs $G_1$ and $G_2$. Then $e\in E(G_1)$ or $e\in E(G_2)$. Without loss of generality, $e\in E(G_1)$. If the join set of the $(\leq 2)$-sum is an edge, then let $vw$ be this edge. Otherwise, let $vw$ be any edge of $G_2$. 

By induction, $G_1$ has a matching $M_1$ with the claimed properties (with respect to the edge $e$), and $G_2$ has a matching $M_2$ with the claimed properties (with respect to the edge $vw$). In particular, every edge in $M_2$ is disjoint from $vw$. Thus $M:=M_1\cup M_2$ is a matching of $G$ (even if $vw\in M_1$). Moreover, every edge in $M$ is disjoint from $e$. We have $|M|= |M_1|+|M_2|\leq \third(\V{G_1}-2)+\third(\V{G_2}-2)= \third(\V{G_1}+\V{G_2}-4)\leq \third(\V{G}-2)$. 

Let $H_i$ be the planar graph obtained by contracting $M_i$ in $G_i$. Let $H$ be the graph obtained by contracting $M$ in $G$. Then $H$ is a $(\leq 2)$-sum of $H_1$ and $H_2$. Thus $H$ is planar.
\end{proof}

\begin{corollary}
\corlabel{PartitionK33}
Every $K_{3,3}$-minor-free graph $G$ has a planar partition with width $2$ and at most $\third(\V{G}-2)$ bags of cardinality $2$. \qed
\end{corollary}

It follows from Euler's Formula and \lemref{CharacterisationK33} that every $K_{3,3}$-minor-free graph $G$ has at most $3\,\V{G}-5$ edges. Thus \corref{PartitionK33} and \lemref{ProduceDrawing} imply the following quantitative version of 
\thmref{RectilinearCrossingK33}.

\begin{corollary}
Every $K_{3,3}$-minor-free graph $G$ has a rectilinear drawing in which each edge crosses at most $2\,\Delta(G)$ other edges. Hence 
\begin{equation*}
\RCR{G}\leq\Delta(G)\,\E{G}\leq \Delta(G)\,(3\,\V{G}-5)\enspace.
\end{equation*}
\qed
\end{corollary}

%%%%%%%%%%%%%%%%%%%%%%%%%%%%%%%%%%%%%%%%%%%%%%%%%%%%%%%%%%%%%%%%%%%%%%%%
\section*{Acknowledgements}
%%%%%%%%%%%%%%%%%%%%%%%%%%%%%%%%%%%%%%%%%%%%%%%%%%%%%%%%%%%%%%%%%%%%%%%%

Thanks to J{\'a}nos Pach and G{\'e}za T{\'o}th for explaining their proofs in reference \cite{PachToth-GD05}. Thanks to Jaroslav Ne{\v{s}}et{\v{r}}il, Patrice Ossona De~Mendez, and Vida Dujmovi\'c for helpful comments.

%%%%%%%%%%%%%%%%%%%%%%%%%%%%%%%%%%%%%%%%%%%%%%%%%%%%%%%%%%%%%%%%%%%%%%%%%%%%%%%%
%\bibliographystyle{nyj}
%\bibliography{myBibliography,myConferences}
%%%%%%%%%%%%%%%%%%%%%%%%%%%%%%%%%%%%%%%%%%%%%%%%%%%%%%%%%%%%%%%%%%%%%%%%%%%%%%%%

\def\Dbar{\leavevmode\lower.6ex\hbox to 0pt{\hskip-.23ex \accent"16\hss}D}

\end{document}